\documentclass[12pt]{amsart}

\usepackage{mathrsfs}
\usepackage{amsfonts}
\usepackage{}
\usepackage{amsmath}
\usepackage{amsthm}
\usepackage{amssymb}
\usepackage{cases}
\usepackage{bbm}
\numberwithin{equation}{section}
\newtheorem{thm}{Theorem}[section]

\newtheorem{lemma}[thm]{Lemma}
\newtheorem{cor}[thm]{Corollary}
\newtheorem{pro}[thm]{Proposition}
\newtheorem{defi}[thm]{Definition}

\newtheorem{prob}[thm]{Problem}

\newcommand{\Rn}{\mathbb{R}^n}
\newcommand{\R}{\mathbb{R}}
\newcommand{\sphere}{S^{n-1}}
\newcommand{\ball}{B^n_2} 
\newcommand{\dom}{\mathrm{dom~} }

\newcommand{\G}{\gamma_{n+1}}
\newcommand{\GD}{\delta\gamma_{n+1} }
\newcommand{\conv}{\mathrm{Conv}}

\newcommand{\DP}{D_{\varphi}}
\newcommand{\DS}{D_{\psi}}
\newcommand{\cn}{c_{n+1}}
\newcommand{\N}{\mathbb{N}}
\usepackage[%
pdfstartview=FitH, %
colorlinks=true, %
linkcolor=blue, %
citecolor=red]{hyperref}

\title [The Gaussian Minkowski problem for epigraphs of convex functions]{The Gaussian Minkowski problem for  epigraphs of convex functions}
\author{$\rm{Xiao\ Li}^1$}
\address{\parbox[l]{1\textwidth}{School of  Mathematical Sciences, Chongqing Normal University,  Chongqing 401331, China}}
\email{lxlixiaolx@163.com}

\author{$\rm{Deping \ Ye}^*$}
\address{\parbox[l]{1\textwidth}{ Department of Mathematics and Statistics, Memorial University of Newfoundland, St.
John's, Newfoundland, A1C 5S7, Canada}}
\email{deping.ye@mun.ca}

\begin{document}

\let\thefootnote\relax
\footnote{Keywords: Convex functions; epigraph; Gaussian moment measure, Minkowski problem, Monge-Amp\`{e}re equation. }

\begin{abstract}  A variational formula is derived by combining the Gaussian volume of the epigraph of a convex function $\varphi$ and the perturbation of $\varphi$ via the infimal convolution. This formula naturally leads to a Borel measure on $\mathbb{R}^n$ and a Borel measure on the unit sphere $S^{n-1}$.  The resulting Borel measure on $\mathbb{R}^n$ will be called the Euclidean Gaussian moment measure of the convex function $\varphi$, and the related   Minkowski-type problem will be studied. In particular, the newly posed Minkowski problem is solved under some mild and natural conditions on the pre-given measure.

\vskip 2mm \noindent {\bf Mathematics Subject Classification (2020):}  26B25 (primary),  52A40, 52A41, 35G20. \end{abstract}

\maketitle
 
 \section{Introduction} Although the term ``Gaussian Minkowski problem" for convex bodies (i.e., compact convex sets in $\Rn$ with nonempty interiors) formally appeared in  \cite{HXZ21} by  Huang, Xi, and Zhao, the problem itself has been posed (albeit implicitly) in \cite{GHWXY19} by Gardner, Hug, Weil, Xing, and Ye. This problem aims to characterize the so-called Gaussian surface area measure of convex bodies. Its normalized version  was first solved in \cite{GHXY20} by Gardner, Hug, Xing, and Ye. In \cite{HXZ21}, Huang, Xi, and Zhao not only provided a solution to the normalized Gaussian Minkowski problem for convex bodies, but more importantly, they provided uniqueness and existence results on the Gaussian Minkowski problem (with no normalization required, which is considerably much more challanging). There is a growing body of work in the Gaussian Minkowski problem and its various extensions  see e.g., \cite{CHLZ23,FHX23,FLX23,KL23, Liu22,SX22,Wang22}. Recently there has been growing attention on the Minkowski-type problems for unbounded closed convex sets. Two typical examples of  unbounded convex sets include the $C$-compatible sets (or $C$-pusedo cones) \cite{AYY24, LYZ24,Schneider2018,Schneider2021,Schneider2024,advSchneider2024, advSchneider2025, S-Zhao,YYZ23}, and the epigraphs of convex functions.  

Our focus in this paper is the epigraphs of convex functions. It is our aim to study the Gaussian Minkowski problem for epigraphs of convex functions, and hence provide a new type of Minkowski problem for convex functions. For convenience, let
$$\conv(\Rn)=\big\{\varphi:\Rn\rightarrow \R\cup\{+\infty\}: \ \varphi\ \ \mathrm{is\ convex, \ lower\mbox- semi}\ \mathrm{continuous}, \varphi\not\equiv +\infty\big\}.$$ By $\dom \varphi$ we mean the  effective domain of $\varphi$ (always convex), i.e., \begin{eqnarray*} 
\dom \varphi =\{x\in\Rn: \varphi(x)<+\infty\}. 
\end{eqnarray*} Thus, $\dom \varphi\neq \emptyset$  if $\varphi\in \conv(\Rn)$.  The epigraph of $\varphi$, denoted by ${\rm epi~} \varphi$, is an unbounded convex set in $\Rn\times\R$  given by: 
\begin{eqnarray*} 
 {\rm epi~} \varphi =\{(x,s)\in\Rn\times\R:\varphi(x)\leq s\}.
\end{eqnarray*} If $\varphi\in \conv(\Rn)$, then ${\rm epi~} \varphi$ is a closed subset in $\R^{n+1}$.

Geometric invariants on epigraphs of convex functions $\varphi\in \conv(\Rn)$ often lead to the functionals on $\varphi$. To see this, following the work \cite{Uli} by Ulivelli, we 
consider a measure $\varpi$ on $\mathbb{R}^{n+1}$ such that
\begin{align*} 
d\varpi(x, s)=\omega (x)\eta(s)dxds,\ \ x\in\Rn \ \ \mathrm{and} \ \  s\in\R,
\end{align*}
where $\omega$ and $\eta$ are nonnegative functions on $\Rn$ and $\R$, respectively.  For $\varphi\in\conv(\Rn)$, let 
\begin{align} 
\varpi(\varphi):=\varpi({\rm epi~} \varphi)=\int_{{\rm epi~} \varphi}d\varpi(x,s)=\int_{\DP }\omega (x)\int_{\varphi(x)}^{+
\infty}\eta(s)ds dx, \label{def-geom-epi-1}
\end{align} where $\DP=\overline{\dom \varphi}$ is the closure of $\dom \varphi$. Some special cases are listed. If $\eta(s)=e^{-s}$, then \eqref{def-geom-epi-1} reduces to the $\omega$-Orlicz moment $\widetilde V_{\omega}(e^{-\varphi})$ defined in \cite{FYZZ23}: \begin{align*} 
\varpi(\varphi)=\int_{\DP }e^{-\varphi(x)}\omega (x)dx,
\end{align*} which includes the total mass (if $\omega(x)\equiv 1$)  and the $(q-n)$-th moment \cite{HLXZJDG} (if $\omega(x)=|x|^{q-n}$ with $|x|$ the Euclidean norm of $x\in \Rn$)  of the log-concave function $e^{-\varphi}$ as its special cases.  If $\omega (x)\equiv 1$ and $\eta(s)=(1-\alpha s)^{\frac{1}{\alpha}-1}\ (-\frac{1}{n}<\alpha<0)$, then \eqref{def-geom-epi-1}  becomes the total mass of the $\alpha$-concave function $(1-\alpha \varphi(x))^{\frac{1}{\alpha}}$ (see e.g., \cite{LNY25,Rot13}) formulated as follows: 
\begin{align*} 
\varpi(\varphi)=\int_{\DP }(1-\alpha \varphi(x))^{\frac{1}{\alpha}}dx. 
\end{align*} Note that, for a convex body $K$, if $\varphi={\rm I}^{\infty}_K$  (taking values $0$ and $+\infty$ on $K$ and outside of $K$, respectively),  by choosing different $\omega$ and $\eta$, $\varpi({\rm I}^{\infty}_K)$ recovers many known geometric invariants on convex bodies, including volume (the total mass of $e^{-{\rm I}^{\infty}_K}$), the $q$-th  dual quermassintegral of $K$ in \cite{Lut75} and the general dual Orlicz quermassintegral of $K$ in \cite{XY-2020-IUMJ, ZXY-2018-JGA}.  In particular, when  $\omega(x)=e^{-\frac{|x|^2}{2}}$, one gets the Gaussian volume of $K$ (up to a multiplicative constant). The Gaussian Minkowski problem in \cite{GHWXY19, HXZ21} aims to find a convex body $K$, such that, for a pre-given Borel measure $\mu$ on the unit sphere $S^{n-1}$, one has  $S_{\gamma_n, K}=\mu.$  Here $S_{\gamma_n, K}$ is the Gaussian surface area measure derived from the following variational formula \cite{GHWXY19}: for two convex bodies $K$ and $L$ containing the origin $o\in \Rn$ in their interiors, one has \begin{align} \label{gauss-vari-geom} \lim_{t\to0^{+}} \frac{(2\pi)^{-\frac{n}{2}}}{t} \bigg(\int _{K +tL} e^{-\frac{|x|^2}{2}}\,dx - \int _{K} e^{-\frac{|x|^2}{2}}\,dx\bigg) =\int_{\sphere}h_LdS_{\gamma_n,K},\end{align} where $K+tL=\{x+ty: x\in K \ \mathrm{and} \ y\in L\}$ for $t>0$, and for a closed (compact or unbounded) convex set $L_1\subset \Rn$,  $h_{L_1}$ denotes its support function taking the following form:  $$h_{L_1}(y)=\sup_{x\in L_1} \langle x, y\rangle, \ \ \mathrm{for} \ \ y\in \Rn, $$
with $\langle x, y\rangle$ being the inner product of $x, y\in \Rn.$ 

The primary goal of this paper is to deal with a Gaussian Minkowski problem for unbounded closed convex sets. More precisely, we are interested in the variational formula for the Gaussian volume of the epigraph of a convex function $\varphi\in \conv(\Rn)$ and related Minkowski problem for epigraphs (and hence for convex functions). Thus, we extend the Gaussian Minkowski problem for convex bodies to convex functions (or some unbounded closed convex sets). By $\G$ we mean the standard Gaussian measure on $\mathbb{R}^{n+1}$, namely, $$\,d\G(x, s)=\cn e^{-\frac{|x|^2+s^2}{2}} \,dx\,ds=\omega(x)\eta(s)\,dx\,ds$$  with $
\omega (x)=(2\pi)^{-\frac{n}{2}}e^{-\frac{|x|^2}{2}}$, $\eta(s)=(2\pi)^{-\frac{1}{2}}e^{-\frac{s^2}{2}}$, and 
\begin{align*} \cn=(2\pi)^{-\frac{n+1}{2}}.
\end{align*}
In this case, we get the Gaussian volume of the epigraph of $\varphi$ (often abbreviated simply as the Gaussian volume of $\varphi$): 
\begin{align*}
\G(\varphi)=\int_{{\rm epi~} \varphi}d\G=\cn \int_{\DP }e^{-\frac{|x|^2}{2}}\int_{\varphi(x)}^{+
\infty}e^{-\frac{s^2}{2}}ds dx.
\end{align*} Clearly, $\G(\varphi)$ is always finite. Note that  \begin{align} \G({\rm I}^{\infty}_K)=\frac{1}{2} (2\pi)^{-\frac{n}{2}}\int_K e^{-\frac{|x|^2}{2}}\,dx= \frac{1}{2}\gamma_n(K),\label{gauss-vol--1} \end{align} where $\gamma_n(K)$ is the Gaussian volume (or measure) of $K$. Due to the nature of the standard Gaussian measure, $\G(\varphi)$ does not have the
translation-invariance and homogeneity. 

In order to setup the Gaussian Minkowski problem for epigraphs, we shall need to define the natural addition for convex functions, which is analogue to the Minkowski addition of convex bodies. Such an addition is called the infimal  convolution  for convex functions $\varphi, \psi\in\mathrm{Conv(\Rn)}$:  
\begin{align*}
\varphi\Box\psi(x)=\inf_{y\in\Rn}\{\varphi(x-y)+\psi(y)\}\quad \text{for } x\in\Rn.
\end{align*} The  right multiplication scalar   of $\varphi$ is defined as
\begin{align*}
(\varphi t)(x)=t \varphi\left(\frac{x}{t}\right)  \quad\text{for }\ t>0\ \ \text{and } x\in\Rn.
\end{align*} The following variation is defined. 

\begin{defi}\label{first-variation}
Let $\varphi,\psi\in \conv(\Rn)$. Define the first variation of the $\varpi(\cdot)$ of $\varphi$ along $\psi$ by  \begin{align*}
\delta_{\varpi} (\varphi,\psi)=\lim_{t\to0^{+}}\frac{\varpi( \varphi\Box(\psi t))-\varpi(\varphi)}{t},
\end{align*} if the limit exists. In particular,  the  first
variation of the Gaussian volume of $\varphi$ along $\psi$ is defined by
\begin{eqnarray}\label{mixedvolume}
\GD(\varphi,\psi)=\lim_{t\to0^{+}}\frac{\G( \varphi\Box(\psi t))-\G(\varphi)}{t}. 
\end{eqnarray} 
\end{defi} 

Before establishing an explicit integral expression for $\GD(\varphi,\psi)$, we briefly review the literature on the results of the integral expressions of $\delta_{\varpi} (\varphi,\psi)$. When $\omega(x) \equiv 1$ and $\eta(s) = e^{-s}$, Klartag and Milman \cite{KM05}  and Rotem \cite{Rot12} studied the special case where $\varphi(x) = \frac{|x|^2}{2}$. Colesanti and Fragalà \cite{CF13} derived integral expressions for the first variation under certain regularity assumptions on $\varphi$ and $\psi$. By using the (anisotropic) coarea formula, these regularity requirements were later removed by Rotem in \cite{Rot20, Rot22} and hence the integral expression of the  first variation has been extended  to more general convex functions.  When $\eta(s) = e^{-s}$, Huang, Liu, Xi, and Zhao \cite{HLXZJDG} obtained the first variation for the $(q-n)$-th moment (i.e., $\omega(x) = |x|^{q-n}$), while Fang, Ye, Zhang, and Zhao \cite{FYZZ23} proved  the first variation for general $\omega$-Orlicz moments, both under certain growth condition near $x = o$.  The additional growth condition (for the $(q-n)$-th moment) was successfully removed by Ulivelli in \cite{Uli}. An $L_p$ version of the first variation for $p > 1$, following the approach in \cite{CF13}, was established by Fang, Xing, and Ye in \cite{FXY}. The approaches in \cite{FYZZ23, HLXZJDG, Rot20, Rot22,  Uli} heavily rely on the variational formulas in geometric settings. As explained in the recent work by Fang, Ye, and Zhang \cite{FYZ24}, an arguably better approach is via analytic techniques and a more suitable set of conditions to impose on $\varphi$ and $\psi$ is arguably the following: there exist constants $\alpha > 0$ and $\beta \in \mathbb{R}$, satisfying that
\begin{align}\label{FYZcondition}
-\infty < \inf \psi^* \leq \psi^* \leq \alpha \varphi^* + \beta \quad \text{on } \mathbb{R}^n, 
\end{align} where $\varphi^*$ denotes the  Legendre transform  of $\varphi$:  
\begin{align*}
\varphi^*(y)=\sup_{x\in\Rn}\left\{\langle x,y\rangle-\varphi(x)\right\}\quad \text{for } y\in\Rn.
\end{align*} It follows from $\psi^*\leq \alpha\varphi^*+\beta$ that $\DS\subseteq \alpha\DP $, which  resembles the condition $L\subseteq a K=\{ax: x\in K\}$, $a>0$, for convex bodies. On the other hand, the condition $-\infty<\inf\psi^*\leq\psi^*$ is used to ensure that $o\in\DS$, resembling the condition $o\in L$ for convex bodies. Under the conditions \eqref{FYZcondition} and that  $o$ is in the interior of $\DP$, Fang, Ye, and Zhang in \cite{FXY} was able to find an integral expression for the first variation of the Riesz $\alpha$-energy for general convex functions $\varphi$ and $\psi$ without the regularity assumptions, the extra growth condition near $x=o$, and the requirement that the effective domain of $\psi$ is a compact set in $\Rn.$ See \cite{FXY} for more details on how to remove the assumption that $o$ is in the interior of $\DP$. The approaches in \cite{FXY} and the condition \eqref{FYZcondition} was successfully used by Li, Nguyen and Ye \cite{LNY25} to calculate the integral expression of the first order variational formula for $\alpha$-concave functions, i.e., $\delta_{\varpi}(\varphi, \psi)$ for $\omega(x)\equiv 1$  and $\eta(s)=(1-\alpha s)^{\frac{1}{\alpha}-1}$ with $-\frac{1}{n}<\alpha <0.$
 
Back to our setting, in Section \ref{section:3},  we will prove a variational formula for the first variation of Gaussian volume of $\varphi$ along $\psi$. For convenience, let $$\mathcal{L}=\Big\{\varphi\in\conv(\Rn): \ \  \liminf\limits_{|x|\to +\infty}\frac{\varphi(x)}{|x|}>0 \Big\}.$$ For $E\subset \Rn$, denote by ${\rm int} (E)$, $\partial E$, and $\mathcal{H}^{n-1}\big|_E$ the interior, boundary, and $(n-1)$-dimensional Hausdorff measure of $E,$ respectively. For the $(n-1)$-dimensional Hausdorff measure of $E,$ we often write  $\mathcal{H}^{n-1}$ if the set $E$ is clearly identified, and in particular, $\,du=\,d\mathcal{H}^{n-1}\big|_{S^{n-1}}(u)$ is often used for the spherical (Lebesgue) measure on the unit sphere $S^{n-1}.$ The set $\DP$ is a closed convex
set, and hence $\partial \DP$ is a Lipschitz manifold, implying that the Gauss map $\nu_{\DP}$ is defined $\mathcal{H}^{n-1}$-almost everywhere 
 on $\partial \DP$. Note that $\varphi\in \conv(\Rn)$ is differentiable almost everywhere in $ \rm int (\DP)$, and when it is differentiable at $x\in \rm int (\DP)$, we shall use $\nabla\varphi(x)$ to denote the gradient of $\varphi$ at $x$. We are now in the position to state our main result in Section \ref{section:3}.  

\vskip 2mm \noindent {\bf Theorem \ref{final}.} {\em  
Let $\varphi\in \mathcal{L}$ be such that  $o\in \rm int (\DP)$. Suppose that $\psi\in \conv(\Rn)$ is a convex function, such that, there exist constants $\alpha>0$ and $\beta\in \R$ satisfying \eqref{FYZcondition}. Then,
\begin{align}\label{vari}
\GD(\varphi,\psi)&=\cn \int_{\partial\DP }h_{\DS}(\nu_{\DP }(x)) e^{-\frac{|x|^2}{2}}\int_{\varphi(x)}^{+
\infty}e^{-\frac{s^2}{2}}ds\,d \mathcal{H}^{n-1}(x)\nonumber\\
&\quad+\cn \int_{\Rn}\psi^*(\nabla \varphi(x)) e^{-\frac{\varphi(x)^2}{2}}e^{-\frac{|x|^2}{2}}dx.
\end{align}}
We point out that Theorem \ref{final} overlaps \cite[Theorem 3.15]{Uli} if we restrict $\varphi$ and $\psi$ to have compact effective domains. The assumption $o\in \rm int(D_{\varphi})$ cannot be removed because $\G(\varphi)$ is not translation invariant when $\varphi$ is replaced by $\varphi(\cdot+x_0)$. We would like to point out that formula \eqref{vari} exhibits the standard structure of the integral expression for the first variation, and these can be seen from similar results in \cite{CF13, FYZ24,FYZZ23, HLXZJDG, LNY25,  Rot22,  Uli}. It induces one Borel measure on $\Rn$ and one Borel measure on $S^{n-1}$. The first one is the push-forward measure of $\cn e^{-\frac{\varphi(x)^2}{2}}e^{-\frac{|x|^2}{2}}\,dx $ under $\nabla \varphi$, which is called  the Euclidean Gaussian moment measure of $\varphi$: for every Borel subset $\vartheta \subseteq \Rn$, \begin{align} 
\mu_{\gamma_{n}}(\varphi,\vartheta)=\cn \int_{\big\{x\in\Rn:\nabla\varphi(x)\in \vartheta \big\}} e^{-\frac{\varphi(x)^2}{2}}e^{-\frac{|x|^2}{2}}dx. \label{euclid-sph-1}
\end{align} The other Borel measure on $S^{n-1}$ is  the push-forward measure, under $\nu_{\DP }$, of  $$\cn \bigg( e^{-\frac{|x|^2}{2}}\int_{\varphi(x)}^{+
\infty}e^{-\frac{s^2}{2}}dsd\mathcal{H}^{n-1}(x)\bigg)\bigg|_{\partial \DP},$$ which is called the spherical Gaussian measure of $\varphi$: for every Borel subset $\vartheta\subseteq \sphere$,
\begin{align} 
\nu_{\gamma_{n}}(\varphi,\vartheta)=\cn \int_{\big\{x\in{\partial\DP }:\ \nu_{\DP }(x)\in \vartheta\big\}} e^{-\frac{|x|^2}{2}}\int_{\varphi(x)}^{+
\infty}e^{-\frac{s^2}{2}}ds\mathcal{H}^{n-1}(x). \label{gauss-sph-1}
\end{align}  By using \eqref{euclid-sph-1}
 and \eqref{gauss-sph-1}, formula \eqref{vari} can be rewritten as follows:  
\begin{align*}
\GD(\varphi,\psi)=
\int_{\Rn}\psi^*(x)d\mu_{\gamma_{n}}(\varphi,x)+\int_{\sphere}h_{\DS}(u)d\nu_{\gamma_{n}}(\varphi,u).
\end{align*} In view of \eqref{gauss-vol--1}, one sees that \eqref{vari} recovers \eqref{gauss-vari-geom}, by letting $\varphi={\rm I}^{\infty}_K$ and $\psi={\rm I}^{\infty}_L$ with $K$ and $L$ two convex bodies  and assuming that condition \eqref{FYZcondition} holds for $\varphi={\rm I}^{\infty}_K$ and $\psi={\rm I}^{\infty}_L$. 

Our second goal in this paper is to study the following Euclidean Gaussian Minkowski problem  for convex functions.

\vskip 2mm\noindent {\bf Problem \ref{problem} {\bf (The Euclidean Gaussian Minkowski problem  for convex functions)}.} {\em 
Let $\mu$ be a nonzero finite Borel measure on $\Rn$. Find the necessary and/or sufficient conditions on $\mu$, such
that,
$$\mu = \tau \mu_{\gamma_n}(\varphi,\cdot)$$
holds for some convex function $\varphi\in\mathcal{L}$ and some constant $\tau>0$.
} 

Although Problem \ref{problem} is stated for convex functions, as previously mentioned, it can be interpreted as a Gaussian Minkowski problem for a family of unbounded closed convex sets (specifically, the epigraphs of convex functions). This formulation extends the Gaussian Minkowski problem for convex bodies \cite{GHWXY19, GHXY20, HXZ21} to  unbounded settings. Once again, finding solutions to Problem \ref{problem} reduces to solving the following Monge-Ampère type equation for an unknown convex function $\varphi$:
\begin{align}\label{Gauss-M-P-M-A}
g(\nabla\varphi(y)){\rm det}(\nabla^2\varphi(y))=\tau\cn  e^{-\frac{\varphi(y)^2}{2}}e^{-\frac{|y|^2}{2}},
\end{align}
where $\tau$ is a constant, ${\rm det}(\nabla^2\varphi(y))$ denotes the determinant of the Hessian matrix of $\varphi$ at $y$, and $d\mu = g(y)\,dy$ with $g$ a smooth function.

Let us pause here to briefly review the literature regarding the Minkowski-type problems for (log-concave, $\alpha$-concave, or convex) functions. As explained before, if $\,d\varpi(x, s)=e^{-s}\,dx\,ds$, then $\varpi(\varphi)$ reduces to the total mass $\int _{\Rn}e^{-\varphi}\,dx$ of $e^{-\varphi}$. The related  Minkowski problem  was initiated by
Cordero-Erausquin and Klartag \cite{CK15} and independently  by  Colesanti and  Fragal\`a \cite{CF13}. Cordero-Erausquin and Klartag \cite{CK15} also obtained the existence and uniqueness of solutions to the functional Minkowski problem aiming to
characterize the moment measure of $\varphi$ (i.e., the push-forward measure of $e^{-\varphi(x)}dx$ under $\nabla \varphi$).  A continuity result for the moment measures has been provided in \cite{Klar14} by  Klartag. Rotem  in \cite{Rot20} and Fang, Xing and Ye in \cite{FXY} provided solutions to the functional $L_p$
Minkowski problem   for $p \in (0, 1)$ and  for $p > 1$,  respectively. The Minkowski problem raised in \cite{CF13} involves two measures (one on $\Rn$ and one on $S^{n-1}$), and recently a solution to this problem has been provided by Falah and Rotem in \cite{Rot25}. The functional dual Minkowski problem (corresponding $\,d\varpi(x,s)=|x|^{q-n}e^{-s}\,dx\,ds$) has been solved in \cite{HLXZJDG} by 
Huang, Liu, Xi and Zhao. In \cite{FYZZ23}, Fang, Ye, Zhang and Zhao solved the functional dual Orlicz Minkowski problem (corresponding $\,d\varpi(x,s)=\omega(x)e^{-s}\,dx\,ds$). Recently, the Riesz $\alpha$-energy Minkowski problem was posed in  \cite{FYZ24} by Fang, Ye and Zhang who also provided a solutions to this problem. These contributions to the solutions for related Minkowski-type problems are primarily based on variational approaches. In particular, for log-concave functions, the identity
$e^{-(\varphi+\psi)} = e^{-\varphi} e^{-\psi}$  
plays a crucial role in solving these problems. This identity  allows translations of $\varphi$ by a constant $a$ (up or down) to a scaling of $e^{-\varphi}$, namely, 
\begin{equation}\label{add}
e^{-(\varphi+a)} = e^{-a} e^{-\varphi}.  
\end{equation}
As a result, certain functionals on log-concave functions, for instance, the total mass, can be easily computed for $e^{-(\varphi+a)}$ and usually have a formulation analogous to   \eqref{add} (probably involving a different power of $e^{-a}$). This property is particularly useful in the variational analysis of the Minkowski-type problems for log-concave functions. It enables the use of the common lower bound 
\begin{align} \label{comm-low-b-1} \varphi(x)\geq a|x|+b  \ \ \mathrm{for} \ \ x\in \Rn\end{align} with $a>0$ and $b\in \R$, without concern for the possibly negative signs of $a|x|+b$ at specific $x\in \Rn$. Moreover, it also allows the transformation of a constrained optimization problem into an unconstrained one, avoiding the need for Lagrange multipliers, which greatly reduces the complexity of solving the related Minkowski problems  (see the details in \cite{CK15, Rot25,FYZ24,FYZZ23,HLXZJDG,Rot20}). Klartag in \cite{Kla18} studied the Minkowski problems for convex functions related to the $q$-moment measure of a convex function $\varphi$, and the Minkowski problem for $\alpha$-concave functions was recently posed and solved in \cite{LNY25} by Li, Nguyen and Ye. The solution to the Minkowski problem for $\alpha$-concave functions in \cite{LNY25} is based on the technique of optimal mass transport,  building upon earlier works by Santambrogio \cite{San16} and by Huynh and Santambrogio \cite{HS21} which dealt with the Minkowski problems for the moment measure and the $q$-moment measure of convex functions, respectively.   

\vskip 1mm Back to our setting, i.e., $\,d\varpi(x, s)=\cn e^{-\frac{|x|^2+s^2}{2}}\,dx\,ds$, in general, one {\em cannot expect} 
 $$\int_{\varphi+a}^{\infty}e^{-\frac{s^2}{2}}\,ds=b\int_{\varphi}^{\infty}e^{-\frac{s^2}{2}}\,ds$$ for some constant $b>0$ (independent of $\varphi$) and hence the identity like  \eqref{add} fails. As a result, to get a non-negative lower bound of $\varphi$,  \eqref{comm-low-b-1} may fail at specific points or regions. Moreover, the transformation of a constrained optimization problem into an unconstrained one is generally not possible. These bring extra difficulty in solving the Euclidean Gaussian Minkowski problem for convex functions (i.e.,  Problem \ref{problem}). These difficulties will be resolved in Section \ref{3.2}, and the proof requires much more work. More specifically, we will replace  \eqref{comm-low-b-1}  by  \begin{align*}  \varphi(x)\geq \max\big \{a|x|+b, \ 0\big\}  \ \ \mathrm{for} \ \ x\in \Rn\end{align*} with $a>0$ and $b\in \R$, and use the method of  Lagrange multipliers to solve Problem \ref{problem}.    
 
 Our solution to  Problem \ref{problem}  is stated and proved in Theorem \ref{main-existence}. 

\vskip 2mm \noindent {\bf Theorem \ref{main-existence}.} {\em
Let $\mu$ be an even nonzero finite Borel measure on $\Rn$ such that $\mu$ is not concentrated in any lower-dimensional subspaces and the first moment of $\mu$ is finite. Then,  there exists  $\varphi\in\mathcal{L}$ such that
\begin{align*} 
 d\mu  =\frac{|\mu|}{\mu_{\gamma_n}(\varphi, \Rn)}d \mu_{\gamma_n}(\varphi, \cdot),
 \end{align*}  where $|\mu|$ and $\mu_{\gamma_n}(\varphi, \Rn)$ are real numbers given by 
 $$\mu_{\gamma_n}(\varphi, \Rn)=\int_{\Rn} d \mu_{\gamma_n}(\varphi, x) \ \ \mathrm{and} \ \ |\mu|=\int_{\Rn}\,d\mu. $$} 

 Note that, in view of \eqref{Gauss-M-P-M-A}, Theorem \ref{main-existence} provides a weak solution to the corresponding Monge-Amp\'{e}re equation. On the other hand, through \eqref{def-geom-epi-1} and  the relations between epigraph and convex function, Theorem \ref{main-existence}  indeed also solves the Gaussian Minkowski problem for some unbounded closed convex sets, which extends those for convex bodies into unbounded settings.

\section{Preliminaries}\label{pre} We now provide some basic definitions and properties for convex functions which are needed in later context. More details can be found in \cite{Roc70,Roc98}. 

Let $\N$ and $\Rn$ be the set of positive integers and the $n$-dimensional Euclidean space with $n\geq 1$, respectively. Denote $o$ the origin in $\Rn.$ A function  $\varphi: \Rn\rightarrow \R\cup \{+\infty\}$ is \emph{convex} if
$$
\varphi((1-\lambda)x+\lambda y)\leq (1-\lambda)\varphi(x)+\lambda \varphi(y),
$$ for all $x,y\in\Rn$ and for $\lambda\in[0,1].$ For a convex function $\varphi$, its   \emph{effective domain}, denoted by $\dom \varphi$, is defined as
\begin{eqnarray*} 
\dom \varphi =\{x\in\Rn: \varphi(x)<+\infty\}.
\end{eqnarray*} Clearly, $\dom \varphi$ is convex in $\Rn$. If  $\dom \varphi\neq\emptyset$, then the convex function $\varphi$ is said to be  \emph{proper}. Let $\DP=\overline{\dom \varphi}$ is the closure of $\dom \varphi$. Associated with convex function $\varphi$  is its \emph{epigraph} ${\rm epi~} \varphi$, a convex set  in $\Rn\times\R$  taking the following form: 
\begin{align*}
 {\rm epi~} \varphi =\{(x,s)\in\Rn\times\R:\varphi(x)\leq s\}.
\end{align*} The set  ${\rm epi~} \varphi$ is  closed, if $\varphi$ is lower semi-continuous. 

Let $\conv(\Rn)$ denote the set of all proper and lower semi-continuous convex functions $\varphi: \Rn\rightarrow \R\cup \{+\infty\}.$ For $\varphi\in \conv(\Rn)$,   ${\rm epi~} \varphi $ must be an unbounded closed convex set, and $\DP$ is also a closed convex set. For a closed convex set $K\subset \Rn$, its boundary $\partial K$ is a Lipschitz manifold and hence  the Gauss map $\nu_{K}$ is well-defined $\mathcal{H}^{n-1}$-almost everywhere 
on $\partial K.$  Hereafter,  $\mathcal{H}^{n-1}|_E$ denotes the $(n-1)$-dimensional Hausdorff measure of the set $E\subset\Rn$, and we often simply use $\mathcal{H}^{n-1}$ if the set $E$ is clearly identified. 
For a set $E\subset\Rn$, by $\overline E$ and $\mathrm{int}(E)$, we mean the closure  and interior of $E$, respectively.
 Let $\omega_n$ denote the volume of the unit ball $\ball$ and $\sphere$ denote the unit sphere. Associated with a closed convex set $K$ is its support function $h_K: \Rn\rightarrow \R$ given by $$h_{K} (y)=\sup_{x\in K}
\langle y, x\rangle \ \ \mathrm{for} \ \ y \in \Rn,$$ with $\langle x, y\rangle$ being the inner product of of $x$ and $y$. In particular, $\nu_{\DP}$ and $h_{\DP}$ are well-defined, and play essential roles in later context. 

The \emph{Legendre transform} $\varphi^*$ of $\varphi$ serves as a natural duality for  a function (not necessarily a convex function)  $\varphi: \Rn\rightarrow \R\cup \{+\infty\}$. It is a convex function of the following form:  
\begin{eqnarray}\label{Fenchel conjugate}
\varphi^*(y)=\sup_{x\in\Rn}\left\{\langle x,y\rangle-\varphi(x)\right\}\quad \text{for } y\in\Rn.
\end{eqnarray} Some easily established results for the Legendre transform are listed here for readers' convenience. Note that for a convex body $K$, 
\begin{align}\label{support}
({\rm I}^{\infty}_K)^*=h_{K},
\end{align}
where ${\rm I}^{\infty}_K$ takes values $0$ and $+\infty$ on $K$ and outside of $K$, respectively. Let $\varphi$ be a proper convex function. Then   \begin{eqnarray}\label{elementary property of Fenchel conjugate}
\varphi^*(o)=-\inf\varphi,
\end{eqnarray} $\varphi^*(y)>-\infty$ for any $y\in\Rn$, and $\varphi^*$ is lower semi-continuous. Moreover,  $\varphi^{**}\leq \varphi$ with equality if and only if $\varphi$ is convex and lower semi-continuous. It also holds that \begin{equation}\label{moto}\varphi^{*}\leq \psi^* \quad \mathrm{if} \quad \varphi\geq\psi.\end{equation} 
It is well known that a proper convex function is continuous in the interior of its effective domain, but differentiable only almost everywhere. When $\varphi$ is differentiable at  $x \in \dom \varphi$, we shall use $\nabla\varphi(x)$ to denote the gradient of $\varphi$ at $x$. Moreover
\begin{eqnarray}\label{basicequ}
\varphi^*(\nabla \varphi(x))+\varphi(x)= \langle x,\nabla \varphi(x)\rangle 
\end{eqnarray} holds at those $x \in \dom \varphi$  where $\varphi$ is differentiable. 

The {\it infimal  convolution}  $\varphi\Box\psi$ of $\varphi, \psi\in \conv(\Rn)$ is defined by 
 \begin{eqnarray}\label{infimal-convolution}
 \varphi\Box\psi(x)=\inf_{y\in\Rn}\{\varphi(x-y)+\psi(y)\} \quad \text{for } x\in\Rn,
 \end{eqnarray}
and the {\it right multiplication scalar} $\varphi t$ of $\varphi\in \conv(\Rn)$ is defined by
 \begin{eqnarray}\label{right scalar multiplication}
(\varphi t)(x)=t \varphi\left(\frac{x}{t}\right)  \quad\text{for } t>0\ \ \text{and } \ x\in\Rn.
\end{eqnarray}  Clearly, these two operations preserve convexity. It can be checked that  
\begin{align*}
 {\rm dom}(\varphi \Box\psi t)&=\dom  \varphi  +t\dom \psi \quad\mathrm{and}\quad
 {\rm epi} (\varphi \Box\psi t)={\rm epi~}  \varphi  +t{\rm epi~} \psi.
\end{align*}The following properties with respect to the operations hold: for $\alpha>0$ and $\beta\in\R$, 
\begin{align} 
 \label{star}
 (\varphi\Box\psi)^* =\varphi^*+\psi^*\quad\mathrm{and}\quad 
 (\psi \alpha-\beta)^* =\alpha\psi^*+\beta.
\end{align}
From \eqref{Fenchel conjugate}, \eqref{moto}, and \eqref{star}, the condition \eqref{FYZcondition} is equivalent to 
\begin{align*}
 o\in\dom\psi \ \ \mathrm{and}\ \ \psi\geq\varphi\alpha-\beta.
\end{align*}

For $\varphi,\psi\in\mathcal{L}$ with  
$$\mathcal{L}=\Big\{\varphi\in\conv(\Rn): \ \  \liminf\limits_{|x|\to +\infty}\frac{\varphi(x)}{|x|}>0 \Big\},$$ where $|x|$ denotes the Euclidean norm of $x\in\Rn$,     $\varphi\Box(\psi t)\in\mathcal{L}$ and thus, 
\begin{align}\label{infimal-Legendre}
\varphi\Box(\psi t)=\big(\varphi\Box(\psi t)\big)^{**}=\big(\varphi^*+t\psi^*\big)^*.
\end{align}
 For $\varphi\in \mathcal{L}$, the condition  $\liminf\limits_{|x|\to +\infty}\frac{\varphi(x)}{|x|}>0 $ implies that there exist constants $a>0$ and $b\in\R$, such that, 
\begin{align}\label{geq}
\varphi(x)\geq a|x|+b \ \ \mathrm{for} \ \ x\in \Rn.\end{align}
Moreover, \eqref{geq} implies that $\int_{\Rn}e^{-\varphi(x)}dx$ is finite, see  e.g.,  {\cite[Lemma 2.5]{CF13}}.

The following result \cite[Proposition 2.1]{Rot20} plays an important role in the later context. 
\begin{lemma}\label{Berman's formula}
Let $\varphi,g:\Rn\to\R\cup\{+\infty\}$ be lower semi-continuous functions. Assume that $g$ is bounded from below and  $g(o),\varphi(o)<+\infty$. Then
$$
\frac{d}{dt}\Big |_{t=0^+}(\varphi+tg)^*(x)=-g(\nabla\varphi^*(x))
$$
at any point $x\in \Rn$ in which $\varphi^*$ is differentiable.
\end{lemma}
If $\varphi,\psi\in\mathcal{L}$ and $\psi^*\geq\inf\psi^*>-\infty$, by \eqref{elementary property of Fenchel conjugate}, one has $\varphi^*(o),\psi^*(o)<+\infty$ and $\psi^*$ is bounded from below. 
Thus the functions $\varphi^*$ and $\psi^*$ satisfy the conditions in Lemma \ref{Berman's formula}, and hence by \eqref{infimal-Legendre}, one has 
\begin{equation}\label{infimal}
\frac{d}{dt}\Big |_{t=0^+}(\varphi\Box(\psi t))(x)=\frac{d}{dt}\Big |_{t=0^+}\big(\varphi^*+t\psi^*\big)^*(x) =-\psi^*(\nabla\varphi(x))
\end{equation}
at any point $x\in \Rn$ in which $\varphi$ is differentiable.

 Recall that, for $\varphi\in\conv(\Rn)$,
\begin{align} 
\G(\varphi)&=\int_{{\rm epi~} \varphi}d\G=\cn \int_{\DP }e^{-\frac{|x|^2}{2}}\int_{\varphi(x)}^{+
\infty}e^{-\frac{s^2}{2}}dsdx \nonumber \\&=\cn  \int_{\Rn }e^{-\frac{|x|^2}{2}}\int_{\varphi(x)}^{+
\infty}e^{-\frac{s^2}{2}}dsdx, \label{volume1}
\end{align}
where $\DP=\overline{\dom \varphi}$  and  $\cn=(2\pi)^{-\frac{(n+1)}{2}}.$

The following lemma holds. 
\begin{lemma}\label{volumefinite}
Let $\varphi\in\conv(\Rn)$. For any $p>0$, one has, 
$$
\int_{\Rn}|x|^p e^{-\frac{|x|^2}{2}}\int_{\varphi(x)}^{+
\infty}e^{-\frac{s^2}{2}}dsdx\in[0,\infty).
$$
\end{lemma}
\begin{proof} Note that the $p$-th moment of the Gaussian measure is finite, which implies  
\begin{align}\label{length}0<\int_{\Rn}|x|^p e^{-\frac{|x|^2}{2}}\,dx<\infty.\end{align}
On the other hand, for any $x\in \Rn,$ $$0\leq \int_{\varphi(x)}^{+
\infty}e^{-\frac{s^2}{2}}ds \leq \sqrt{2\pi}.$$ These yield that $$0\leq \int_{\Rn}|x|^p e^{-\frac{|x|^2}{2}}\int_{\varphi(x)}^{+
\infty}e^{-\frac{s^2}{2}}dsdx\leq \sqrt{2\pi} \int_{\Rn}|x|^p e^{-\frac{|x|^2}{2}}dx<\infty.$$ This concludes the proof. 
\end{proof}

We shall also need the following result. 
\begin{lemma}\label{volumefinite2}
Let $\varphi\in\conv(\Rn)$. Then
$$
\bigg|\int_{\Rn}\varphi(x)e^{-\frac{|x|^2}{2}}e^{-\frac{\varphi(x)^2}{2}}dx\bigg| \leq 
\int_{\Rn} \big|\varphi(x)\big|e^{-\frac{|x|^2}{2}}e^{-\frac{\varphi(x)^2}{2}}dx<\infty.
$$
\end{lemma}
\begin{proof} It is easily checked that, for $t\geq 0,$ $$te^{-\frac{t^2}{2}}\leq {e}^{-\frac{1}{2}}.$$ By letting $t=|\varphi(x)|$, one gets 
\begin{align*}
\int_{\Rn} \big|\varphi(x)\big|e^{-\frac{|x|^2}{2}}e^{-\frac{\varphi(x)^2}{2}}dx & \leq   {e}^{-\frac{1}{2}} \int_{\Rn}  e^{-\frac{|x|^2}{2}}dx  <\infty. 
\end{align*}  This completes the proof. 
\end{proof}

We now prove the  last result in this section. 

\begin{lemma}\label{volumefinite1}
Let $\varphi\in\mathcal{L}$. Then
$$
\int_{\Rn}|\nabla \varphi(x)| e^{-\frac{\varphi(x)^2}{2}}e^{-\frac{|x|^2}{2}}dx\in[0,\infty).
$$
\end{lemma}
\begin{proof}
It has been proved in \emph{\cite[Lemma 4]{CK15}} that, for $\varphi\in\mathcal{L}$,  
$$\int_{\Rn}|\nabla e^{-\varphi(x)}|dx=\int_{\Rn}|\nabla \varphi(x)| e^{-\varphi(x)}dx\in[0,\infty).$$ This further yields that 
\begin{eqnarray*}\label{measurepro}
\int_{\Rn}|\nabla\varphi(x)|e^{-\frac{\varphi(x)^2}{2}}e^{-\frac{|x|^2}{2}}dx
&\leq& \int_{\Rn}|\nabla\varphi(x)|e^{-\frac{\varphi(x)^2}{2}}dx\nonumber\\
&\leq& e^{\frac{1}{2}}\int_{\Rn}|\nabla\varphi(x)|e^{-\varphi(x)}dx<\infty,
\end{eqnarray*} where we have used the inequality  $\frac{r^2}{2}\geq r-\frac{1}{2}$ for any $r$.
\end{proof}
 
\section{A variational formula for the Gaussian volume of the epigraphs of convex functions} \label{section:3}
 In this section, we will calculate the explicit integral expression for  $\GD(\varphi,\psi)$  defined in \eqref{mixedvolume}: 
for  $\varphi,\psi\in \conv(\Rn)$,
$$
\GD(\varphi,\psi)=\lim_{t\to0^{+}}\frac{\G( \varphi\Box(\psi t))-\G(\varphi)}{t}.
$$  
Let us first prove the following property for $\GD(\varphi,\psi)$:
\begin{pro}\label{cor}
Let  $\varphi,\psi\in \mathcal{L}$ be such that $ \psi^*\geq\inf\psi^*>-\infty$  and   $\GD(\varphi,\psi)$ exist. Assume that, for some $\alpha>0$ and $\beta\in\R$, $\widetilde\psi=\psi \alpha-\beta$  satisfies  that $$\lim_{t\rightarrow0^+}\varphi\Box(\widetilde\psi t)=\varphi.$$ Then, the following holds:  $$\GD(\varphi,\widetilde\psi)=\alpha\GD(\varphi,\psi)+
\beta\cn \int_{\Rn}e^{-\frac{|x|^2}{2}}e^{-\frac{\varphi(x)^2}{2}}dx.$$
\end{pro}

\begin{proof}
Set $\overline\varphi_t=\varphi\Box(\widetilde\psi t)$. From \eqref{infimal-convolution} and \eqref{right scalar multiplication}, one has
\begin{equation}\label{double}
\overline\varphi_t=\varphi\Box(\psi(\alpha t))-t\beta,
\end{equation} where $\psi(\alpha t)$ is the right multiplication of $\psi$ and $\alpha t$. 
Since $\lim_{t\rightarrow0^+}\varphi\Box(\widetilde\psi t)=\varphi$, then $\lim_{t\rightarrow 0^+}\varphi\Box(\psi(\alpha t))=\varphi$. Based on \eqref{volume1}, we can rewrite 
\begin{align}\GD(\varphi,\widetilde\psi)&=\lim_{t\to0^{+}} \frac{\G(\overline\varphi_t)-\G( \varphi)}{t}  =\cn (B_1+B_2), \label{formula-B1}
\end{align}
where $B_1$ and $B_2$ are given by: 
$$B_1=\lim_{t\to0^{+}} \frac{1}{t}\int_{\Rn}e^{-\frac{|x|^2}{2}} 
\bigg(
 {\int_{\overline\varphi_t(x)}^{+
\infty}e^{-\frac{s^2}{2}}ds-\int_{\varphi\Box(\psi(\alpha t))(x)}^{+
\infty}e^{-\frac{s^2}{2}}ds}\bigg)dx,$$
$$B_2=\lim_{t\to0^{+}} \frac{1}{t}  \int_{\Rn} e^{-\frac{|x|^2}{2}}\bigg(\int_{\varphi\Box(\psi(\alpha t))(x)}^{+
\infty}e^{-\frac{s^2}{2}}ds-\int_{\varphi(x)}^{+
\infty}e^{-\frac{s^2}{2}}ds\bigg)dx.$$
From \eqref{double}, we can get $$\frac{1}{t}\left| \int_{\overline\varphi_t(x)}^{+\infty}e^{-\frac{s^2}{2}}ds-\int_{\varphi\Box(\psi(\alpha t))(x)}^{+
\infty}e^{-\frac{s^2}{2}}ds \right|\leq  |\beta|.$$
Together with \eqref{star} and \eqref{infimal}, the dominated convergence theorem gives that
\begin{align}
B_1&=\int_{\Rn}e^{-\frac{|x|^2}{2}} \lim_{t\to0^{+}} \bigg[\frac{1}{t} 
\bigg(
 {\int_{\overline\varphi_t(x)}^{+
\infty}e^{-\frac{s^2}{2}}ds-\int_{\varphi\Box(\psi(\alpha t))(x)}^{+
\infty}e^{-\frac{s^2}{2}}ds}\bigg)\bigg]dx \nonumber\\ &=\int_{\Rn}(-\alpha\psi^*(\nabla \varphi(x))) e^{-\frac{|x|^2}{2}}e^{-\frac{\varphi(x)^2}{2}}dx-\int_{\Rn}(-\alpha\psi^*-\beta)(\nabla \varphi(x))e^{-\frac{|x|^2}{2}}e^{-\frac{\varphi(x)^2}{2}}dx \nonumber\\
&=\beta\int_{\Rn}e^{-\frac{|x|^2}{2}}e^{-\frac{\varphi(x)^2}{2}}dx. \label{formula-B2}
\end{align}
According to  \eqref{mixedvolume} for $\GD(\varphi,\psi)$,  one has \begin{align} 
\cn  B_2 &= \alpha \lim_{t\to0^{+}}\frac{\G( \varphi\Box(\psi (\alpha t)))-\G(\varphi)}{\alpha t}  \nonumber \\&= \alpha \lim_{\tau \to0^{+}}\frac{\G( \varphi\Box(\psi \tau))-\G(\varphi)}{\tau}  \nonumber \\  &= \alpha\GD(\varphi,\psi), \label{formula-B3} \end{align} where we have used the substitution $\tau=\alpha t.$ The desired formula follows directly from \eqref{formula-B1}, \eqref{formula-B2} and \eqref{formula-B3}.
\end{proof}

When $\widetilde\psi=\varphi\alpha-\beta$, it follows from \eqref{infimal-convolution} and \eqref{right scalar multiplication} that, for any $x\in \Rn$, \begin{align}\label{beishu}
\varphi\Box(\widetilde\psi t)(x)=\varphi\Box((\varphi\alpha-\beta) t)(x)=(1+\alpha t)\varphi\bigg(\frac{x}{1+\alpha t}\bigg)-t\beta.\end{align}
If $o\in \rm int (\DP)$,  one has $\varphi^*(y)\geq-\varphi(o)>-\infty$ for any $y\in\Rn$ due to \eqref{Fenchel conjugate}. Moreover, it follows from the lower semi-continuity of $\varphi$, \cite[Lemma 1.6.11]{Schneider2014} and \eqref{beishu} that 
\begin{align}\label{conti-varphi}
\lim_{t\rightarrow0^+}\varphi\Box((\varphi\alpha-\beta) t)=\varphi.
\end{align} 
Hence we can immediately get the following result.
\begin{cor}\label{cor1} Let $\varphi\in \mathcal{L}$ be such that  $o\in \rm int (\DP)$ and $\GD(\varphi,\varphi)$ exist. Then, for $\alpha>0$ and $\beta\in\R$, one has 
$$\GD(\varphi,\varphi \alpha-\beta)=\alpha\GD(\varphi,\varphi)+
\beta\cn \int_{\Rn}e^{-\frac{|x|^2}{2}}e^{-\frac{\varphi(x)^2}{2}}dx.$$
\end{cor}

Subsequently, we will calculate $\GD(\varphi,\psi)$ following the proofs of the first order variational formula for the Riesz $\alpha$-energy \cite{FYZ24} and the total mass of $\alpha$-concave functions \cite{LNY25}. Firstly,  we calculate  $\GD(\varphi,\varphi).$
\begin{lemma}\label{variation formula}
Let $\varphi \in \conv(\Rn)$. Then
\begin{eqnarray*}
\GD(\varphi,\varphi)\!=\!n\G( \varphi)\!-\!\cn \bigg(\int_{\Rn}\!|x|^2e^{-\frac{|x|^2}{2}}\!\!\int_{\varphi(x)}^{+
\infty}\!\!\!\!e^{-\frac{s^2}{2}}dsdx\!+\!\int_{\Rn}\!\!\varphi(x)e^{-\frac{|x|^2}{2}}e^{-\frac{\varphi(x)^2}{2}}dx\bigg).
\end{eqnarray*}
In particular, $\GD(\varphi,\varphi)$ is finite.
\end{lemma}

\begin{proof} Note that $(\varphi\Box(\varphi t))(x)=(1+t)\varphi(\frac{x}{1+t})$. It follows from \eqref{volume1} that 
\begin{align*}\label{}
\G(\varphi\Box(\varphi t))
&=\cn \int_{\Rn}e^{-\frac{|x|^2}{2}}\int_{(1+t)\varphi(\frac{x}{1+t})}^{+
\infty}e^{-\frac{s^2}{2}}dsdx\nonumber\\
&=\cn (1+t)^n\int_{\Rn}e^{-\frac{|(1+t)z|^2}{2}}\int_{(1+t)\varphi(z)}^{+
\infty}e^{-\frac{s^2}{2}}dsdz.
\end{align*}
where we used the substitution $x=(1+t)z$.
This further implies that 
\begin{align}\label{zong}
  \cn^{-1} \GD(\varphi,\varphi) &=      \lim_{t\to0^{+}} \frac{\G(\varphi\Box(\varphi t))-\G( \varphi)}{\cn t}\notag\\
&=  \lim_{t\to0^{+}} \frac{1}{t}\bigg(\!(1+t)^n \!\! \int_{\Rn}\!\!\!e^{-\frac{|(1+t)x|^2}{2}} \! \! \!\int_{(1+t)\varphi(x)}^{+
\infty} \! \! \!\!\!\!\! e^{-\frac{s^2}{2}}dsdx-\!\int_{\Rn}\!\! e^{-\frac{|x|^2}{2}}\!\! \int_{\varphi(x)}^{+
\infty}\!\!\! e^{-\frac{s^2}{2}}dsdx\! \bigg)\notag\\
&=A_1+A_2+A_3,
\end{align}
where $A_1, A_2$ and $A_3$ are given by  \begin{align*} 
A_1&=\lim_{t\to0^{+}}\frac{(1+t)^n-1}{t}\int_{\Rn}e^{-\frac{|(1+t)x|^2}{2}}\int_{(1+t)\varphi(x)}^{+
\infty}e^{-\frac{s^2}{2}}dsdx,\\ A_2 &=\lim_{t\to0^{+}}\int_{\Rn}\frac{e^{-\frac{|(1+t)x|^2}{2}}-e^{-\frac{|x|^2}{2}}}{t}\int_{(1+t)\varphi(x)}^{+
\infty}e^{-\frac{s^2}{2}}dsdx,\\ A_3&=\lim_{t\to0^{+}}\int_{\Rn}e^{-\frac{|x|^2}{2}} \frac{1}{t} \bigg({\int_{(1+t)\varphi(x)}^{+
\infty}e^{-\frac{s^2}{2}}ds-\int_{\varphi(x)}^{+
\infty}e^{-\frac{s^2}{2}}ds}\bigg)dx. \end{align*}
Note that, for any $t>0$ and $x\in\Rn$,
$$e^{-\frac{|(1+t)x|^2}{2}}\int_{(1+t)\varphi(x)}^{+
\infty}e^{-\frac{s^2}{2}}ds\leq (2\pi)^{\frac{1}{2}}e^{-\frac{|x|^2}{2}}.$$ It follows from  the dominated convergence theorem  that
\begin{align}\label{A-1}
A_1&=\lim_{t\to0^{+}}\frac{(1+t)^n-1}{t}\int_{\Rn}\lim_{t\to0^{+}}e^{-\frac{|(1+t)x|^2}{2}}\int_{(1+t)\varphi(x)}^{+
\infty}e^{-\frac{s^2}{2}}dsdx\notag\\
&=n\int_{\Rn}e^{-\frac{|x|^2}{2}}\int_{\varphi(x)}^{+
\infty}e^{-\frac{s^2}{2}}dsdx.
\end{align}

Next we compute $A_2.$ By the mean value theorem, for $0\leq t\leq1$ and  $x\in\Rn$, there exists $s_x\in(0,t)$, such that, 
\begin{align*} 
0  \leq e^{-\frac{|x|^2}{2}}-e^{-\frac{|(1+t)x|^2}{2}}  =t|x|^2  (1+s_x) e^{-\frac{(1+s_x)^2|x|^2}{2}}     \leq2t|x|^2e^{-\frac{|x|^2}{2}}.\end{align*}

Since $\int_{(1+t)\varphi(x)}^{+
\infty}e^{-\frac{s^2}{2}}ds\leq (2\pi)^{\frac{1}{2}}$, we can get
$$\left|\frac{e^{-\frac{|(1+t)x|^2}{2}}-e^{-\frac{|x|^2}{2}}}{t}\int_{(1+t)\varphi(x)}^{+
\infty}e^{-\frac{s^2}{2}}ds\right|\leq 2^{\frac{3}{2}}\pi^{\frac{1}{2}}|x|^2e^{-\frac{|x|^2}{2}}.$$
Together with \eqref{length}, the dominated convergence theorem yields that
\begin{align}\label{A-2}
A_2&=\int_{\Rn}\lim_{t\to0^{+}}\bigg[\frac{e^{-\frac{|(1+t)x|^2}{2}}-e^{-\frac{|x|^2}{2}}}{t}\int_{(1+t)\varphi(x)}^{+
\infty}e^{-\frac{s^2}{2}}ds\bigg]dx\notag\\
&=\int_{\Rn}\lim_{t\to0^{+}}\frac{e^{-\frac{|(1+t)x|^2}{2}}-e^{-\frac{|x|^2}{2}}}{t}\lim_{t\to0^{+}}\int_{(1+t)\varphi(x)}^{+
\infty}e^{-\frac{s^2}{2}}dsdx\notag\\
&=-\int_{\Rn}|x|^2e^{-\frac{|x|^2}{2}}\int_{\varphi(x)}^{+
\infty}e^{-\frac{s^2}{2}}dsdx.
\end{align}

Finally we calculate $A_3$. 
From   $te^{-\frac{t^2}{2}}\leq e^{-\frac{1}{2}}$ for $t\geq 0,$ we can get
\begin{align*}
 \frac{1}{t}\bigg| \int_{(1+t)\varphi(x)}^{\varphi(x)}e^{-\frac{s^2}{2}}ds \bigg|\leq \left|e^{-\frac{\varphi(x)^2}{2}}\frac{t\varphi(x)}{t}\right|\leq e^{-\frac{1}{2}}.
\end{align*}
It follows from the dominated convergence theorem that
\begin{align}\label{A-3}
A_3&=\int_{\Rn}e^{-\frac{|x|^2}{2}}\lim_{t\to0^{+}} \frac{1}{t} \bigg({\int_{(1+t)\varphi(x)}^{+
\infty}e^{-\frac{s^2}{2}}ds-\int_{\varphi(x)}^{+
\infty}e^{-\frac{s^2}{2}}ds}\bigg)dx\notag\\&=-\int_{\Rn}\varphi(x)e^{-\frac{|x|^2}{2}}e^{-\frac{\varphi(x)^2}{2}}dx.\end{align}
The conclusion follows from \eqref{zong}, \eqref{A-1}, \eqref{A-2} and \eqref{A-3}.
\end{proof}

Next we give an integral formula of $\GD(\varphi,\varphi)$.
\begin{pro}\label{integral formula} 
Let $\varphi\in \mathcal{L}$ be such that  $o\in \rm int (\DP)$. Then
\begin{align*} 
\GD(\varphi,\varphi)&=\cn \int_{\partial\DP }\langle x,\nu_{\DP }(x)\rangle e^{-\frac{|x|^2}{2}}\int_{\varphi(x)}^{+
\infty}e^{-\frac{s^2}{2}}ds \,d \mathcal{H}^{n-1}(x)\\
&\quad+\cn \int_{\Rn}\varphi^*(\nabla \varphi(x)) e^{-\frac{\varphi(x)^2}{2}}e^{-\frac{|x|^2}{2}}dx.
\end{align*}
\end{pro}
\begin{proof}
From $o\in \rm int(\DP)$ and the convexity of $\varphi$, it holds that 
\begin{align}\label{compact}
\langle x,\nabla\varphi(x)\rangle\geq \varphi(x)-\varphi(o)\geq\inf \varphi-\varphi(o)>-\infty,
\end{align}
at any point $x\in \Rn$ in which $\varphi$ is differentiable.

Let $B_2^n(R)$ denote the ball with radial $R$ centered at the origin and $\mathrm {div}$ be the divergence operator. It follows from  \eqref{compact} and the monotone convergence theorem (applied to the nonnegative function $\langle x, \nabla \varphi\rangle+\varphi(o)-\inf \varphi\geq 0$) that
\begin{align} 
\int_{\Rn}\langle x,\nabla\varphi(x)\rangle e^{-\frac{\varphi(x)^2}{2}}e^{-\frac{|x|^2}{2}}dx
&=-\int_{\Rn}\Big\langle x,\nabla\int_{\varphi(x)}^{+
\infty}e^{-\frac{s^2}{2}}ds\Big\rangle e^{-\frac{|x|^2}{2}}dx \nonumber \\
&=-\lim_{R\rightarrow \infty}\int_{\DP \cap B_2^n(R)}\Big\langle xe^{-\frac{|x|^2}{2}},\nabla\int_{\varphi(x)}^{+
\infty}e^{-\frac{s^2}{2}}ds\Big\rangle dx\nonumber \\ &=\lim_{R\rightarrow \infty}\int_{\DP \cap B_2^n(R)}{\mathrm {div}}\Big( x e^{-\frac{|x|^2}{2}}\Big)\int_{\varphi(x)}^{+
\infty}e^{-\frac{s^2}{2}}dsdx \nonumber \\
&\quad-\lim_{R\rightarrow \infty}\int_{\DP \cap B_2^n(R)}\!\!\!  {\mathrm {div}}\Big( x e^{-\frac{|x|^2}{2}}\int_{\varphi(x)}^{+
\infty} \!\!\!  e^{-\frac{s^2}{2}}ds\Big)dx. \label{3.15-1}
\end{align}
It can be checked that, for any $x\in \Rn,$ $${\mathrm {div}}\Big( x e^{-\frac{|x|^2}{2}}\Big)=ne^{-\frac{|x|^2}{2}}-|x|^2e^{-\frac{|x|^2}{2}}\geq ne^{-\frac{|x|^2}{2}}-\frac{4}{e}e^{-\frac{|x|^2}{4}}.$$ 
Applying the monotone convergence theorem to the following nonnegative function $$\Big({\mathrm {div}}\Big( x e^{-\frac{|x|^2}{2}}\Big)-ne^{-\frac{|x|^2}{2}}+\frac{4}{e}e^{-\frac{|x|^2}{4}}\Big) \int_{\varphi(x)}^{+
\infty}e^{-\frac{s^2}{2}}ds\geq 0,$$ one can deduce that 
 \begin{align}   \lim_{R\rightarrow \infty}\!\!\int_{\DP \cap B_2^n(R)}\!\!\!\!\! {\mathrm {div}}\Big( x e^{-\frac{|x|^2}{2}}\Big)\!\!\int_{\varphi(x)}^{+
\infty}\!\! \!\! e^{-\frac{s^2}{2}}dsdx & = \int_{\Rn}{\mathrm {div}}\Big( x e^{-\frac{|x|^2}{2}}\Big)\int_{\varphi(x)}^{+
\infty}e^{-\frac{s^2}{2}}dsdx \nonumber \\&=\int_{\Rn}(n-|x|^2)e^{-\frac{|x|^2}{2}}\int_{\varphi(x)}^{+
\infty}e^{-\frac{s^2}{2}}dsdx. \label{3.15-2} 
\end{align}
We now claim that
\begin{align}\label{3.6}
\!\!\!\!\!\! \lim_{R\rightarrow \infty}\! \int_{\DP \cap B_2^n(R)}\!\!\!\!\!\!\!\! {\mathrm {div}}\Big( x e^{-\frac{|x|^2}{2}}\!\!\!\int_{\varphi(x)}^{+
\infty}\!\!\!\!\! e^{-\frac{s^2}{2}}ds\Big)dx
\!=\!\!\int_{\partial\DP }\!\!\!\!\langle x,\nu_{\DP }(x)\rangle e^{-\frac{|x|^2}{2}}\!\!\!\int_{\varphi(x)}^{+
\infty}\!\!\!e^{-\frac{s^2}{2}}dsd\mathcal{H}^{n-1}(x).
\end{align}
 The divergence theorem can be applied to get
\begin{align}\label{divide-DP}
& \lim_{R\rightarrow \infty}\int_{\DP \cap B_2^n(R)}{\mathrm {div}}\Big( x e^{-\frac{|x|^2}{2}}\int_{\varphi(x)}^{+
\infty}e^{-\frac{s^2}{2}}ds\Big)dx\notag\\  
&\quad =\lim_{R\rightarrow \infty}\int_{\partial(\DP \cap B_2^n(R))}\langle x,\nu_{\DP \cap B_2^n(R)}(x)\rangle e^{-\frac{|x|^2}{2}}\int_{\varphi(x)}^{+
\infty}e^{-\frac{s^2}{2}}dsd\mathcal{H}^{n-1}(x)\notag\\
& \quad =\lim_{R\rightarrow \infty}\int_{\Xi_1(R)}\langle x,\nu_{B_2^n(R)}(x)\rangle e^{-\frac{|x|^2}{2}}\int_{\varphi(x)}^{+
\infty}e^{-\frac{s^2}{2}}dsd\mathcal{H}^{n-1}(x)\notag\\
&\quad\quad +\lim_{R\rightarrow \infty}\int_{\Xi_2(R)}\langle x,\nu_{\DP }(x)\rangle e^{-\frac{|x|^2}{2}}\int_{\varphi(x)}^{+
\infty}e^{-\frac{s^2}{2}}dsd\mathcal{H}^{n-1}(x),
\end{align}
where $\Xi_1(R)=\partial(\DP \cap B_2^n(R))\cap\partial(B_2^n(R))$ and $\Xi_2(R)=\partial(\DP \cap B_2^n(R))\cap\partial\DP $.

Direct computation gives that  
\begin{align*}
&0\leq\lim_{R\rightarrow \infty}\int_{\Xi_1(R)}\langle x,\nu_{B_2^n(R)}(x)\rangle e^{-\frac{|x|^2}{2}}\int_{\varphi(x)}^{+
\infty}e^{-\frac{s^2}{2}}dsd\mathcal{H}^{n-1}(x)\notag\\
&\leq(2\pi)^{\frac{1}{2}}\lim_{R\rightarrow \infty}\int_{\sphere}R^n e^{-\frac{R^2}{2}}du\notag\\
&=(2\pi)^{\frac{1}{2}}n\omega_n\lim_{R\rightarrow \infty}R^n e^{-\frac{R^2}{2}} =0.
\end{align*} 
Consequently, it follows that 
\begin{align}\label{3.7}
&\lim_{R\rightarrow \infty}\int_{\Xi_1(R)}\langle x,\nu_{B_2^n(R)}(x)\rangle e^{-\frac{|x|^2}{2}}\int_{\varphi(x)}^{+
\infty}e^{-\frac{s^2}{2}}dsd\mathcal{H}^{n-1}(x)=0.
\end{align} 
On the other hand, by $o\in \rm int(\DP)$, the Cauchy-Schwarz inequality,   $e^{-\frac{s^2}{2}}\leq e^{\frac{1}{2}}e^{-s}$ for all $s\in\R$, and $ e^{\frac{1}{2}}se^{-\frac{s^2}{2}}\leq 1$ for $s>0$, 
one has,  for $x\in\partial \DP$, 
\begin{align}\label{partialdom}
0 \leq\langle x,\nu_{\DP }(x)\rangle e^{-\frac{|x|^2}{2}}\int_{\varphi(x)}^{+
\infty} e^{-\frac{s^2}{2}}ds  \leq |x| e^{-\frac{|x|^2}{2}}  e^{\frac{1}{2}} \int_{\varphi(x)}^{+
\infty} e^{-s}ds   \leq e^{-\varphi(x)}.
\end{align}
It follows from \eqref{geq} and \eqref{partialdom} that
\begin{align*}
&0\leq\lim_{R\rightarrow \infty}\int_{\partial\DP \setminus\Xi_2(R)}\langle x,\nu_{\DP }(x)\rangle e^{-\frac{|x|^2}{2}}\int_{\varphi(x)}^{+
\infty}e^{-\frac{s^2}{2}}dsd\mathcal{H}^{n-1}(x)\\
&\leq\lim_{R\rightarrow \infty}\int_{\partial\DP \setminus\Xi_2(R)}e^{-\varphi(x)}d\mathcal{H}^{n-1}(x)\\
&\leq e^{-\frac{b}{2}}\lim_{R\rightarrow \infty} e^{-\frac{aR}{2}}\int_{\partial\DP }e^{-\frac{\varphi(x)}{2}}d\mathcal{H}^{n-1}(x)\\
&=0,
\end{align*}
for some $a>0$ and $b\in\R$, where the last equality follows from the fact proved in \cite[Proposition 1.6]{Rot22}, but applied to $\frac{\varphi}{2}$, that $\int_{\partial{\DP}}e^{-\frac{\varphi(x)}{2}}d\mathcal{H}^{n-1}(x)$ is finite. This further implies the following identity: 
\begin{align*}
\lim_{R\rightarrow \infty}\!\int_{\Xi_2(R)}\!\!\!\!\!\langle x,\nu_{\DP }(x)\rangle e^{-\frac{|x|^2}{2}}\!\!\!\int_{\varphi(x)}^{+
\infty}\!\!\!\!e^{-\frac{s^2}{2}}dsd\mathcal{H}^{n-1}(x)\!=\!\!\int_{\partial\DP }\!\!\!\!\!\langle x,\nu_{\DP }(x)\rangle e^{-\frac{|x|^2}{2}}\!\!\!\int_{\varphi(x)}^{+
\infty}\!\!\!\!e^{-\frac{s^2}{2}}dsd\mathcal{H}^{n-1}(x).
\end{align*} This together with   \eqref{divide-DP} and \eqref{3.7} yields the claim \eqref{3.6}. 

Combining \eqref{3.15-1}, \eqref{3.15-2}, and \eqref{3.6}, it follows that  
 \begin{align*}
& n\int_{\Rn}\!\! e^{-\frac{|x|^2}{2}}\!\! \int_{\varphi(x)}^{+\infty}\!\! \!\!\! e^{-\frac{s^2}{2}}dsdx\!-\!\! \int_{\Rn}\!\!|x|^2e^{-\frac{|x|^2}{2}}\!\int_{\varphi(x)}^{+
\infty}\!\!\!\!\! e^{-\frac{s^2}{2}}dsdx  \nonumber \\ & \quad \!=\! \int_{\Rn}\!\!\langle x,\nabla\varphi(x)\rangle e^{-\frac{\varphi(x)^2}{2}}e^{-\frac{|x|^2}{2}}dx  +\!\int_{\partial\DP }\!\langle x,\nu_{\DP }(x)\rangle e^{-\frac{|x|^2}{2}}\!\!\int_{\varphi(x)}^{+
\infty}\!\!\!\!\! e^{-\frac{s^2}{2}}dsd\mathcal{H}^{n-1}(x).
\end{align*}
Together with Lemma \ref{variation formula} and \eqref{basicequ},
one has
\begin{align*}
\GD(\varphi,\varphi)
&=\cn \int_{\partial\DP }\langle x,\nu_{\DP }(x)\rangle e^{-\frac{|x|^2}{2}}\int_{\varphi(x)}^{+
\infty}e^{-\frac{s^2}{2}}dsd\mathcal{H}^{n-1}(x)\\
&\quad +\cn \int_{\Rn}\langle x,\nabla\varphi(x)\rangle e^{-\frac{\varphi(x)^2}{2}}e^{-\frac{|x|^2}{2}}dx
-\cn \int_{\Rn}\varphi(x)e^{-\frac{|x|^2}{2}}e^{-\frac{\varphi(x)^2}{2}}dx\\
&=\cn \int_{\partial\DP }\langle x,\nu_{\DP }(x)\rangle e^{-\frac{|x|^2}{2}}\int_{\varphi(x)}^{+
\infty}e^{-\frac{s^2}{2}}dsd\mathcal{H}^{n-1}(x)\\
&\quad +\cn \int_{\Rn}\varphi^*(\nabla \varphi(x)) e^{-\frac{\varphi(x)^2}{2}}e^{-\frac{|x|^2}{2}}dx.
\end{align*}
This completes the proof. 
\end{proof}

The following lemma is needed to establish the explicit integral expression of $\GD(\varphi,\psi)$.
\begin{lemma}\label{lem6}
Let $\varphi\in \mathcal{L}$ be such that  $o\in \rm int (\DP)$ and $\widehat\varphi_t=\varphi\Box((\varphi\alpha-\beta) t)$ for some $\alpha>0$ and $\beta\in\R$. For $u\in\sphere$, one has \begin{align*}
        \lim_{t\to 0^{+}}\int_{\sphere}E_{t}(u)du = \int_{\sphere}\lim_{t\to 0^{+}} E_{t}(u)du<\infty,
    \end{align*} where $E_{t}: S^{n-1}\rightarrow \R$ is defined by 
    \begin{align}\label{et}E_{t}(u) := 
\frac{1}{t}\int_{0}^{+\infty}e^{-\frac{r^2}{2}}r^{n-1}\int_{\widehat\varphi_t(ru)}^{\varphi(ru)}e^{-\frac{s^2}{2}}dsdr.\end{align} 
     
\end{lemma}
\begin{proof}
    By repeating the proof of Lemma \ref{variation formula} and by Corollary \ref{cor1}, one has
\begin{align*}
\lim_{t\to 0^{+}}E_{t}(u)&= \lim_{t\to 0^{+}}\int_{0}^{+\infty}e^{-\frac{r^2}{2}}r^{n-1}\frac{1}{t}\int_{\widehat\varphi_t(ru)}^{\varphi(ru)}e^{-\frac{s^2}{2}}dsdr\nonumber\\
    &=\alpha\bigg( n\!\!\int_{0}^{+\infty}\!\!\!e^{-\frac{r^2}{2}}r^{n-1}\!\!\int_{\varphi(ru)}^{+
\infty}\!\!e^{-\frac{s^2}{2}}dsdr\! -\!\!\int_{0}^{+\infty}\!\!e^{-\frac{r^2}{2}}r^{n+1}\!\!\int_{\varphi(ru)}^{+\infty}\!\!e^{-\frac{s^2}{2}}dsdr\bigg)\nonumber\\
&\quad-\alpha\!\!\int_{0}^{+\infty}\!\!\varphi(ru)e^{-\frac{r^2}{2}}e^{-\frac{\varphi(ru)^2}{2}}r^{n-1}dr+\beta\!\!\int_{0}^{+\infty}\!\!e^{-\frac{r^2}{2}}e^{-\frac{\varphi(ru)^2}{2}}r^{n-1}dr.
\end{align*}
It follows from the polar coordinate formula and Lemma \ref{variation formula} that
\begin{align}\label{Et}
\!\!  \int_{\sphere}\lim_{t\to 0^{+}}E_{t}(u)du &=\beta\int_{\Rn}e^{-\frac{|x|^2}{2}}e^{-\frac{\varphi(x)^2}{2}}dx-\alpha\int_{\Rn}\varphi(x)e^{-\frac{|x|^2}{2}}e^{-\frac{\varphi(x)^2}{2}}dx\notag\\& \quad +\alpha\bigg(n\int_{\Rn }e^{-\frac{|x|^2}{2}}\int_{\varphi(x)}^{+
\infty}\!\!\!\! e^{-\frac{s^2}{2}}dsdx-\!\! \int_{\Rn}|x|^2e^{-\frac{|x|^2}{2}}\int_{\varphi(x)}^{+
\infty}\!\!\!\! e^{-\frac{s^2}{2}}dsdx\!\bigg)\notag\\
&=\cn^{-1}\bigg(\alpha\GD(\varphi,\varphi)+
\beta\cn \int_{\Rn}e^{-\frac{|x|^2}{2}}e^{-\frac{\varphi(x)^2}{2}}dx\bigg).
\end{align} 
    Due to the polar coordinate formula   and Corollary \ref{cor1}, we have
    \begin{align}\label{3.91}
        \lim_{t\to 0^{+}}\int_{\sphere}E_{t}(u)du &=\lim_{t\to 0^{+}}\int_{\sphere}\int_{0}^{+\infty}e^{-\frac{r^2}{2}}r^{n-1}\int_{\widehat\varphi_t(ru)}^{\varphi(ru)}e^{-\frac{s^2}{2}}dsdrdu\nonumber\\
&=\!\lim_{t\to 0^{+}}\int_{\Rn}e^{-\frac{|x|^2}{2}}\int_{\widehat\varphi_t(x)}^{\varphi(x)}e^{-\frac{s^2}{2}}dsdx\nonumber\\    
 &=\!\cn^{-1}\lim_{t\to 0^{+}}\frac{\G(\widehat\varphi_t)-\G(\varphi)}{t}\nonumber\\
        &=\cn^{-1}\bigg(\alpha\GD(\varphi,\varphi)+
\beta\cn \int_{\Rn}e^{-\frac{|x|^2}{2}}e^{-\frac{\varphi(x)^2}{2}}dx\bigg).
        \end{align}
 The conclusion follows from \eqref{Et} and \eqref{3.91}.
\end{proof}

If $o\in {\rm int( \DP)}$,  for $u\in\sphere$, one can define $\rho_{\DP }: S^{n-1}\rightarrow [0, +\infty]$, the radial function of $\DP $ (not necessarily compact), by
$$\rho_{\DP }(u) := \sup\{t >0 : tu \in\DP \}.$$
For $u\in\sphere$, set
\begin{align} A_t(u)=\frac{1}{t}\int_{\rho_{\DP }(u)}^{\rho_{D_{\widehat\varphi_t}}(u)}e^{-\frac{r^2}{2}}r^{n-1}\int_{\widehat\varphi_t(ru)}^{+
\infty}e^{-\frac{s^2}{2}}dsdr, \nonumber \\ B_t(u)=\frac{1}{t}\int_{0}^{\rho_{\DP }(u)}e^{-\frac{r^2}{2}}r^{n-1}\int_{\widehat\varphi_t(ru)}^{\varphi(ru)}e^{-\frac{s^2}{2}}dsdr. \label{def-B-t} \end{align}
Then, $E_t(u)$ can be rewritten as
\begin{align*}E_t(u)&=A_t(u)+B_t(u).\end{align*}
Therefore by \eqref{Et}, one has 
\begin{align}\label{BT}
\int_{\sphere}\lim_{t\to 0^{+}}E_t(u)du&=\int_{\sphere}\lim_{t\to 0^{+}}\big(A_t(u)+B_t(u)\big)du \notag\\
&=\cn^{-1}\bigg(\alpha\GD(\varphi,\varphi)+
\beta\cn \int_{\Rn}e^{-\frac{|x|^2}{2}}e^{-\frac{\varphi(x)^2}{2}}dx\bigg).\end{align}

\begin{lemma}\label{equality formula}
Let $\varphi\in \mathcal{L}$ be such that  $o\in \rm int (\DP)$ and $\widehat\varphi_t=\varphi\Box((\varphi\alpha-\beta) t)$ for some $\alpha>0$ and $\beta\in\R$.  Then, for almost every $u\in\sphere$, one has
\begin{align*}
\lim_{t\to0^{+}}\frac{1}{t}\int_{0}^{\rho_{\DP }(u)}\!\!\!e^{-\frac{r^2}{2}}r^{n-1}\int_{\widehat\varphi_t(ru)}^{\varphi(ru)}e^{-\frac{s^2}{2}}dsdr
=\!\int_{0}^{\rho_{\DP }(u)}\!\!e^{-\frac{r^2}{2}}r^{n-1}\lim_{t\to0^{+}}\frac{1}{t}\int_{\widehat\varphi_t(ru)}^{\varphi(ru)}e^{-\frac{s^2}{2}}dsdr.
\end{align*}
\end{lemma}

\begin{proof}
Let $\Omega_{\varphi}=\{u\in\sphere:\rho_{\DP }(u)<+\infty\}$. If $u\notin \Omega_{\varphi}$, one has $\rho_{\DP }(u)=+\infty$, and then   $\rho_{D_{\widehat\varphi_t }} (u)=+\infty$ by $D_{\widehat\varphi_t }=(1+\alpha t)\DP$. Therefore for $u\notin \Omega_{\varphi}$,
\begin{align}\label{infinite}
A_t(u)=\frac{1}{t}\int_{\rho_{\DP }(u)}^{\rho_{D_{\widehat\varphi_t }}(u)}e^{-\frac{r^2}{2}}r^{n-1}\int_{\widehat\varphi_t(ru)}^{+
\infty}e^{-\frac{s^2}{2}}dsdr=0.\end{align}

Subsequently we consider $u\in \Omega_{\varphi}$, that is  $\rho_{\DP }(u)<+\infty$. 
By variable change $r=\tau \rho_{D_{\varphi}}(u)$ with $\tau  \in [1, 1+\alpha t]$ and the mean
value theorem for the definite integrals, there exists $ \tau_0(t,u)\in(1,1+\alpha t)$, such that, 
\begin{align}\label{atu}
\lim_{t\to0^{+}}A_t(u)&=\lim_{t\to0^{+}}\frac{1}{t}\left(\int_{\rho_{\DP }(u)}^{\rho_{D_{\widehat\varphi_t }}(u)}e^{-\frac{r^2}{2}}r^{n-1}\int_{\widehat\varphi_t(ru)}^{+
\infty}e^{-\frac{s^2}{2}}dsdr\right)\notag\\
&=\rho^n_{\DP }(u)\lim_{t\to0^{+}}\frac{1}{t}\left(\int_{1}^{1+\alpha t}e^{-\frac{\left(\tau\rho_{\DP }(u)\right)^2}{2}}\tau^{n-1}\int_{\widehat\varphi_t(\tau\rho_{\DP }(u)u)}^{+
\infty}e^{-\frac{s^2}{2}}dsd\tau\right)\notag\\
&=\alpha\rho^n_{\DP }(u)\lim_{t\to0^{+}}e^{-\frac{\left(\tau_0(t,u)\rho_{\DP }(u)\right)^2}{2}}\tau_0(t,u)^{n-1}\int_{\widehat\varphi_t(\tau_0(t,u)\rho_{\DP }(u)u)}^{+
\infty}e^{-\frac{s^2}{2}}ds.
\end{align}

Note that $\lim_{t\to0^{+}}\tau_0(t,u)=1^+$. According to \eqref{conti-varphi}, for $u\in\sphere$ and $0<r<\rho_{\DP }(u)$, one has (see a detailed argument on page 22 in \cite{FYZ24}),
\begin{equation}\label{double1continuity}
\lim_{t\to0^{+}}\widehat\varphi_t(ru)=\varphi(ru)\ \ \ \mathrm{and}\ \  \  \lim_{t\to0^{+}}\widehat\varphi_t(\tau_0(t,u)\rho_{\DP}(u)u)=\varphi(\rho_{\DP}(u)u).
\end{equation} 
Together with \eqref{atu}, we get, for $u\in\Omega_{\varphi}$, \begin{align}\label{A1}
\lim_{t\to0^{+}}A_t(u)
&=\alpha\rho^n_{\DP }(u)e^{-\frac{\left(\rho_{\DP }(u)\right)^2}{2}}\int_{\varphi(\rho_{\DP }(u)u)}^{+
\infty}e^{-\frac{s^2}{2}}ds.
\end{align} 
It follows from \eqref{infinite}, \eqref{A1}, and the variable change $x=\rho_{D_{\varphi}}(u)u$  that 
\begin{align*}
\int_{\sphere}\lim_{t\to0^{+}}A_t(u)du
&=\alpha\int_{\Omega_{\varphi}}\rho^n_{\DP }(u)e^{-\frac{\left(\rho_{\DP }(u)\right)^2}{2}}\int_{\varphi(\rho_{\DP }(u)u)}^{+
\infty}e^{-\frac{s^2}{2}}dsdu\notag\\
&=\alpha\int_{\partial\DP }\langle x,\nu_{\DP }(x)\rangle e^{-\frac{|x|^2}{2}}\int_{\varphi(x)}^{+
\infty}e^{-\frac{s^2}{2}}ds\,d\mathcal{H}^{n-1}(x).
\end{align*} Together with  Proposition \ref{integral formula}  and \eqref{BT}, one has 
\begin{align}\label{keylemmma}
\!\!\int_{\sphere}\lim_{t\to 0^{+}} B_t(u)du &=\int_{\sphere}\lim_{t\to 0^{+}} \frac{1}{t}\int_{0}^{\rho_{\DP }(u)}e^{-\frac{r^2}{2}}r^{n-1}\int_{\widehat\varphi_t(ru)}^{\varphi(ru)}e^{-\frac{s^2}{2}}dsdrdu  \notag\\
& =
\cn^{-1}\bigg(\alpha\GD(\varphi,\varphi)\!+\!
\beta\cn \!\! \int_{\Rn}\!\!e^{-\frac{|x|^2}{2}}e^{-\frac{\varphi(x)^2}{2}}dx\bigg)\!-\!\!\int_{\sphere}\lim_{t\to 0^{+}} \! A_t(u)du\notag\\
& =\alpha\int_{\Rn}e^{-\frac{|x|^2}{2}}\varphi^*(\nabla \varphi(x)) e^{-\frac{\varphi(x)^2}{2}}dx
+\beta\int_{\Rn}e^{-\frac{|x|^2}{2}}e^{-\frac{\varphi(x)^2}{2}}dx.
\end{align}

On the other hand, by \eqref{star}, \eqref{infimal} and \eqref{double1continuity},  for almost all $u\in\sphere$ and $0<r<\rho_{\DP }(u)$, we have 
\begin{align*}
&\int_{0}^{\rho_{\DP }(u)}\!\!e^{-\frac{r^2}{2}}r^{n-1}\lim_{t\to0^{+}}\frac{1}{t}\int_{\widehat\varphi_t(ru)}^{\varphi(ru)}e^{-\frac{s^2}{2}}dsdr\notag\\&\quad=\alpha\!\!\int_{0}^{\rho_{\DP }(u)}\!\!e^{-\frac{r^2}{2}}r^{n-1}\varphi^*(\nabla \varphi(ru)) e^{-\frac{\varphi(ru)^2}{2}}dr
\!+\!\beta \!\int_{0}^{\rho_{\DP }(u)}\!\!\!e^{-\frac{r^2}{2}}r^{n-1}e^{-\frac{\varphi(ru)^2}{2}}dr.
\end{align*}
Together with the polar coordinate formula, one has
\begin{align}\label{A2}
&\int_{\sphere}\int_{0}^{\rho_{\DP }(u)}e^{-\frac{r^2}{2}}r^{n-1}\lim_{t\to0^{+}}\frac{1}{t}\int_{\widehat\varphi_t(ru)}^{\varphi(ru)}e^{-\frac{s^2}{2}}dsdrdu\notag\\&\quad=\alpha\int_{\Rn}e^{-\frac{|x|^2}{2}}\varphi^*(\nabla \varphi(x)) e^{-\frac{\varphi(x)^2}{2}}dx
+\beta\int_{\Rn}e^{-\frac{|x|^2}{2}}e^{-\frac{\varphi(x)^2}{2}}dx.
\end{align}
The desired equality follows from \eqref{keylemmma} and \eqref{A2}. \end{proof}

We are now in the position to prove our main theorem in this section, following a similar approach to those of \cite[Theorem 1.4]{FYZ24} and \cite[Theorem 3.10]{LNY25}.  
\begin{thm}\label{final}
Let $\varphi\in \mathcal{L}$ be such that  $o\in \rm int (\DP)$. Suppose that $\psi\in \conv(\Rn)$ is a convex function, such that, there exist constants $\alpha>0$ and $\beta\in \R$ satisfying \eqref{FYZcondition}. Then,
\begin{align}\label{vari--1}
\GD(\varphi,\psi)&=\cn \int_{\partial\DP }h_{\DS}(\nu_{\DP }(x)) e^{-\frac{|x|^2}{2}}\int_{\varphi(x)}^{+
\infty}e^{-\frac{s^2}{2}}ds\,d \mathcal{H}^{n-1}(x)\nonumber\\
&\quad+\cn \int_{\Rn}\psi^*(\nabla \varphi(x)) e^{-\frac{\varphi(x)^2}{2}}e^{-\frac{|x|^2}{2}}dx.
\end{align}
\end{thm}

\begin{proof}
Write $\varphi_t=\varphi\Box(\psi t)$ and $\widehat\varphi_t=\varphi\Box((\varphi\alpha-\beta) t)$. First, we assume that $\inf\psi^*\geq0$. It follows from \eqref{FYZcondition} that, for $t\geq 0,$ $$\varphi^*\leq\varphi^*+t\psi^*\leq (1+\alpha t)\varphi^*+\beta t.$$
 By \eqref{moto}, \eqref{star} and \eqref{infimal-Legendre}, one has
\begin{align}\label{baohan}\widehat\varphi_t\leq\varphi_t\leq\varphi \quad \rm{and}\quad \DP \subseteq D_{\varphi_t}\subseteq D_{\widehat\varphi_t}.\end{align}
This together with \eqref{double1continuity} implies that, for $x\in \rm int (\DP),$
\begin{align}\label{continuity}
\lim_{t\to0^{+}}\varphi_t(x)=\varphi(x).
\end{align} Combining \eqref{et} and \eqref{baohan}, one has, for $u\in\sphere$, 
$$0\leq\frac{1}{t}\int_{0}^{+\infty}e^{-\frac{r^2}{2}}r^{n-1}\int_{\varphi_t(ru)}^{\varphi(ru)}e^{-\frac{s^2}{2}}dsdr\leq E_t(u).$$
It follows from  Lemma \ref{lem6} and the general dominated convergence theorem  that 
\begin{align}
\GD(\varphi,\psi)&=\lim_{t\to0^{+}} \frac{\G(\varphi_t)-\G( \varphi)}{t}\notag\\
&=\cn \lim_{t\to 0^{+}}\frac{1}{t}\int_{\sphere}\!\int_{0}^{+\infty}e^{-\frac{r^2}{2}}r^{n-1}\int_{\varphi_t(ru)}^{\varphi(ru)}e^{-\frac{s^2}{2}}dsdrdu\notag\\
&=\cn \int_{\sphere}\lim_{t\to 0^{+}}\frac{1}{t}\!\int_{0}^{+\infty}e^{-\frac{r^2}{2}}r^{n-1}\int_{\varphi_t(ru)}^{\varphi(ru)}e^{-\frac{s^2}{2}}dsdrdu\notag\\
&=\cn \int_{\sphere}\lim_{t\to0^{+}}\big(C_t(u)+D_t(u)\big)du,\label{jiao}
\end{align}
where $C_t$ and $D_t$ are given by 
\begin{align*}
    C_t(u)&=\frac{1}{t}\int_{\rho_{\DP }(u)}^{\rho_{D_{\varphi_t}}(u)}e^{-\frac{r^2}{2}}r^{n-1}\int_{\varphi_t(ru)}^{+
\infty}e^{-\frac{s^2}{2}}dsdr,\\ D_t(u)&=\frac{1}{t}\int_{0}^{\rho_{\DP }(u)}e^{-\frac{r^2}{2}}r^{n-1}\int_{\varphi_t(ru)}^{\varphi(ru)}e^{-\frac{s^2}{2}}dsdr. 
\end{align*}

First we compute $\lim_{t\to0^{+}}C_t(u)$. By $D_{\varphi_t}=\DP +t\DS$ and \cite[Lemma 5.3]{FYZ24}, one has, for $u\in\Omega_{\varphi}$,  $$\lim_{t\to0^{+}}\frac{\rho_{D_{\varphi_t}}(u)-\rho_{\DP }(u)}{t}=\frac{h_{\DS}(\nu_{\DP }(\rho_{\DP }(u)u))}{\langle u,\nu_{\DP }(\rho_{\DP }(u)u)\rangle}.$$ The mean
value theorem for the definite integrals  and \eqref{continuity} yield that, for $u\in\Omega_{\varphi}$,
\begin{align}\label{Ct}
\lim_{t\to0^{+}}C_t(u)&=\lim_{t\to0^{+}}\frac{1}{t}\int_{\rho_{\DP }(u)}^{\rho_{D_{\varphi_t}}(u)}e^{-\frac{r^2}{2}}r^{n-1}\int_{\varphi_t(ru)}^{+
\infty}e^{-\frac{s^2}{2}}dsdr\notag\\
&=\lim_{t\to0^{+}}\bigg(\frac{\rho_{D_{\varphi_t}}(u)-\rho_{\DP }(u)}{t}
 \tau(t,u)^{n-1}e^{-\frac{\left(\tau(t,u)\right)^2}{2}}\int_{\varphi_t(\tau(t,u)u)}^{+
\infty}e^{-\frac{s^2}{2}}ds\bigg)\notag\\ &=\bigg(\lim_{t\to0^{+}}\frac{\rho_{D_{\varphi_t}}(u)-\rho_{\DP }(u)}{t}\bigg)
\bigg(\lim_{t\to0^{+}}\tau(t,u)^{n-1}e^{-\frac{\left(\tau(t,u)\right)^2}{2}}\int_{\varphi_t(\tau(t,u)u)}^{+
\infty}e^{-\frac{s^2}{2}}ds\bigg)\notag\\
&=\frac{h_{\DS}(\nu_{\DP }(\rho_{\DP }(u)u))}{\langle u,\nu_{\DP }(\rho_{\DP }(u)u)\rangle}
\rho^{n-1}_{\DP }(u)e^{-\frac{\left(\rho_{\DP }(u)\right)^2}{2}}\int_{\varphi(\rho_{\DP }(u)u)}^{+
\infty}e^{-\frac{s^2}{2}}ds,
\end{align}
where $\tau(t,u)\in(\rho_{\DP }(u),\rho_{D_{\varphi_t}}(u))$ satisfies $\lim_{t\to0^{+}}\tau(t,u)=\rho_{\DP }(u)$. Similar to \eqref{infinite}, one has, for $u\notin\Omega_\varphi$,
\begin{align}\label{ct2}\lim_{t\to0^{+}}C_t(u)=\lim_{t\to0^{+}}\frac{1}{t}\int_{\rho_{\DP }(u)}^{\rho_{D_{\varphi_t}}(u)}e^{-\frac{r^2}{2}}r^{n-1}\int_{\varphi_t(ru)}^{+
\infty}e^{-\frac{s^2}{2}}dsdr=0.\end{align}

Second, let us deal with $\lim_{t\to0^{+}}D_t(u)$.
It follows from \eqref{baohan} that, for almost all $u\in\sphere$ and $0<r<\rho_{\DP }(u)$,   
$0\leq D_t(u)\leq B_t(u)$ with $B_t$ given in \eqref{def-B-t}. That is, 
$$0\leq\int_{0}^{\rho_{\DP }(u)}e^{-\frac{r^2}{2}}r^{n-1}\int_{\varphi_t(ru)}^{\varphi(ru)}e^{-\frac{s^2}{2}}dsdr\leq \int_{0}^{\rho_{\DP }(u)}e^{-\frac{r^2}{2}}r^{n-1}\int_{\widehat\varphi_t(ru)}^{\varphi(ru)}e^{-\frac{s^2}{2}}dsdr.$$
Together with Lemma \ref{equality formula}, the general dominated convergence theorem implies that
\begin{align*}
\lim_{t\to0^{+}}D_t(u)&=\lim_{t\to0^{+}}\frac{1}{t}\int_{0}^{\rho_{\DP }(u)}e^{-\frac{r^2}{2}}r^{n-1}\int_{\varphi_t(ru)}^{\varphi(ru)}e^{-\frac{s^2}{2}}dsdr\\
&=\int_{0}^{\rho_{\DP }(u)}e^{-\frac{r^2}{2}}r^{n-1}\lim_{t\to0^{+}}\frac{1}{t}\int_{\varphi_t(ru)}^{\varphi(ru)}e^{-\frac{s^2}{2}}dsdr.
\end{align*}
It follows from \eqref{infimal} and \eqref{continuity} that, for almost all $u\in\sphere$ and $0<r<\rho_{\DP }(u)$,
\begin{align*}
\lim_{t\to0^{+}}\frac{1}{t}\int_{\varphi_t(ru)}^{\varphi(ru)}e^{-\frac{s^2}{2}}ds=\psi^*(\nabla \varphi(ru)) e^{-\frac{\varphi(ru)^2}{2}}.
\end{align*}
Hence, for almost all $u\in\sphere$, one has \begin{align}\label{Dt}
\lim_{t\to0^{+}}D_t(u)=\int_{0}^{\rho_{\DP }(u)}e^{-\frac{r^2}{2}}r^{n-1}\psi^*(\nabla \varphi(ru)) e^{-\frac{\varphi(ru)^2}{2}}dr.
\end{align}

Combining \eqref{jiao}, \eqref{Ct}, \eqref{ct2}, \eqref{Dt}, and the polar coordinate formula, one gets 
\begin{align}
\GD(\varphi,\psi)
&=\cn \int_{\Omega_\varphi}\frac{h_{\DS}(\nu_{\DP }(\rho_{\DP }(u)u))}{\langle u,\nu_{\DP }(\rho_{\DP }(u)u)\rangle}
\rho^{n-1}_{\DP }(u)e^{-\frac{\left(\rho_{\DP }(u)\right)^2}{2}}\int_{\varphi(\rho_{\DP }(u)u)}^{+
\infty}e^{-\frac{s^2}{2}}dsdu \nonumber  \\
&\quad+\cn \int_{\sphere}\int_{0}^{\rho_{\DP }(u)}e^{-\frac{r^2}{2}}r^{n-1}\psi^*(\nabla \varphi(ru)) e^{-\frac{\varphi(ru)^2}{2}}drdu \nonumber  \\
&=\cn \int_{\partial\DP }h_{\DS}(\nu_{\DP }(x)) e^{-\frac{|x|^2}{2}}\int_{\varphi(x)}^{+
\infty}e^{-\frac{s^2}{2}}dsd\mathcal{H}^{n-1}(x)\nonumber \\
&\quad+\cn \int_{\Rn}\psi^*(\nabla \varphi(x)) e^{-\frac{\varphi(x)^2}{2}}e^{-\frac{|x|^2}{2}}dx.\label{desired-positive-1}
\end{align} This shows \eqref{vari--1} when $\inf \psi^*\geq 0.$

Finally, we deal with the case when $\inf\psi^*<0$. Set $\widetilde \psi^*=\psi^*-\inf\psi^*$, which yields $\widetilde \psi =\psi +\inf\psi^*$ due to  \eqref{star}. Then, ${\DS}=D_{\widetilde\psi}$ and
\begin{align*}
0\leq\widetilde \psi^*\leq \alpha\varphi^*+\big(\beta-\inf\psi^*\big).
\end{align*} Similar to \eqref{continuity}, one has 
$$
\lim_{t\to0^{+}}\varphi\Box(\widetilde\psi t)=\varphi.
$$ It follows from Proposition \ref{cor} that 
\begin{align*}
\GD(\varphi,\widetilde\psi)   
 =\GD(\varphi,\psi)-
\inf\psi^*\cn \int_{\Rn}e^{-\frac{|x|^2}{2}}e^{-\frac{\varphi(x)^2}{2}}dx. 
\end{align*} Applying \eqref{desired-positive-1} to $\widetilde \psi$ (satisfying $\widetilde \psi^*\geq 0$), one has  
\begin{align*}
\GD(\varphi,\psi)   
&=\GD(\varphi,\widetilde\psi) +
\inf\psi^*\cn \int_{\Rn}e^{-\frac{|x|^2}{2}}e^{-\frac{\varphi(x)^2}{2}}dx\\ 
&=\cn \left(\int_{\Rn}\widetilde \psi^*(\nabla \varphi(x)) e^{-\frac{\varphi(x)^2}{2}}e^{-\frac{|x|^2}{2}}dx+\inf\psi^*\int_{\Rn} e^{-\frac{\varphi(x)^2}{2}}e^{-\frac{|x|^2}{2}}dx\right)\\
&\quad+\cn \int_{\partial\DP }h_{D_{\widetilde \psi}}(\nu_{\DP }(x)) e^{-\frac{|x|^2}{2}}\int_{\varphi(x)}^{+
\infty}e^{-\frac{s^2}{2}}dsd\mathcal{H}^{n-1}(x)\\&= \cn  \int_{\Rn}  \psi^*(\nabla \varphi(x)) e^{-\frac{\varphi(x)^2}{2}}e^{-\frac{|x|^2}{2}}dx  \\
&\quad+\cn \int_{\partial\DP }h_{\DS}(\nu_{\DP }(x)) e^{-\frac{|x|^2}{2}}\int_{\varphi(x)}^{+
\infty}e^{-\frac{s^2}{2}}dsd\mathcal{H}^{n-1}(x).   
\end{align*} Consequently, the desired result holds, and this completes the proof.  
\end{proof}

Theorem \ref{final} motivates two Borel measures as defined below. 
\begin{defi}\label{def-two-measures}
Let $\varphi\in \mathcal{L}$ be a convex function. 

\vskip 2mm \noindent i) The Euclidean Gaussian moment measure $\mu_{\gamma_n}(\varphi,\cdot)$ of $\varphi$ is a Borel measure on $ \Rn$ defined as follows: for every Borel subset $\vartheta\subseteq \Rn$,
\begin{eqnarray*}
\mu_{\gamma_n}(\varphi,\vartheta)=\cn \int_{\big\{x\in\Rn:\ \nabla\varphi(x)\in \vartheta\big\}} e^{-\frac{\varphi(x)^2}{2}}e^{-\frac{|x|^2}{2}}dx,
\end{eqnarray*}
where $\nabla\varphi$ is the gradient of $\varphi$, i.e.,  $\mu_{\gamma_n}(\varphi,\cdot)$ is the push-forward measure of $\cn e^{-\frac{\varphi(x)^2}{2}}e^{-\frac{|x|^2}{2}}\,dx $ under the map $\nabla \varphi$.

\vskip 2mm \noindent ii)  The spherical Gaussian moment measure $\nu_{\gamma_n}(\varphi,\cdot)$ of $\varphi$ is a Borel measure on $ \sphere$ defined as follows: for every Borel subset $\vartheta\subseteq \sphere$,
\begin{eqnarray*}
\nu_{\gamma_n}(\varphi,\vartheta)=\cn \int_{\big\{x\in{\partial\DP }:\ \nu_{\DP }(x)\in \vartheta 
\big\}} e^{-\frac{|x|^2}{2}}\int_{\varphi(x)}^{+
\infty}e^{-\frac{s^2}{2}}dsd\mathcal{H}^{n-1}(x).
\end{eqnarray*}
 where $\nu_{\DP}$ is the Gauss map of $\partial \DP$. That is,  $\nu_{\gamma_n}(\varphi,\cdot)$  is the push-forward measure (on the unit sphere $S^{n-1}$) of  $\cn e^{-\frac{|x|^2}{2}}\int_{\varphi(x)}^{+
\infty}e^{-\frac{s^2}{2}}dsd\mathcal{H}^{n-1}(x)|_{\partial \DP}$  under the map $\nu_{\DP }$.
\end{defi}

 Using the above notations, one can rewrite \eqref{vari--1} as
\begin{eqnarray*}
\GD(\varphi,\psi)=
\int_{\Rn}\psi^*(x)d\mu_{\gamma_n}(\varphi,x)+\int_{\sphere}h_{D_\psi}(u)d\nu_{\gamma_n}(\varphi,u).
\end{eqnarray*}
If $K,L$ are two convex bodies with $o\in {\rm int }(K)$ and $o\in L$, then ${\rm I}^{\infty}_K$ and ${\rm I}^{\infty}_L$ satisfy the condition \eqref{FYZcondition}. That is, from \eqref{support}, the condition \eqref{FYZcondition} is equivalent to the following fact:  \begin{align*}
-\infty< \inf h_L\leq h_L \leq \alpha h_K \ \  \mathrm{on}\ \ \Rn,
\end{align*} for some constant $\alpha>0$. By Theorem \ref{final},  one has
\begin{align*}
\GD({\rm I}^{\infty}_K,{\rm I}^{\infty}_L)
&=\lim_{t\to0^{+}}\frac{\G({\rm I}^{\infty}_{K+L} )-\G({\rm I}^{\infty}_K)}{t}\\
&=\frac{1}{2}\lim_{t\to0^{+}}\frac{\gamma_{n}(K+L)-\gamma_{n}(K)}{t}\\
&=\frac{1}{2}(2\pi)^{-\frac{n}{2}}\int_{\partial{K}}h_{L}(\nu_{K}(x)) e^{-\frac{|x|^2}{2}}d\mathcal{H}^{n-1}(x). 
\end{align*} Consequently,  Theorem \ref{final} recovers the variational formula \eqref{gauss-vari-geom} of Gaussian
volume $\gamma_{n}(K)$ obtained in \cite{GHWXY19, GHXY20,  HXZ21}.

Definition \ref{def-two-measures} motivates the following  Minkowski-type problem:

\begin{prob}  Let $\mu$ and $\nu$ be finite Borel measures on $\Rn$ and $\sphere$, respectively. What are the  necessary and/or sufficient conditions on $\mu$ and $\nu$ so that there exist some convex functions $\varphi\in\mathcal{L}$ and constants $\tau_1, \tau_2$  satisfying 
$$\mu =\tau_1 \mu_{\gamma_n}(\varphi,\cdot)\ \ \ \mathrm{and} \ \ \ \nu=\tau_2\nu_{\gamma_n}(\varphi,\cdot). $$
\end{prob} 

In Section \ref{3.2}, we shall concentrate on the special case when $\nu$ is a zero measure. That is, we aim to solve the following Minkowski problem regarding the Euclidean Gaussian moment measure $\mu_{\gamma_n}(\varphi,\cdot)$. 
\begin{prob}[{\bf The Euclidean Gaussian Minkowski problem  for convex functions}]  \label{problem} Let $\mu$ be a nonzero finite Borel measure on $\Rn$. Find the necessary and/or sufficient conditions on $\mu$, such that 
$$\mu = \tau \mu_{\gamma_n}(\varphi,\cdot)$$
holds for some convex function $\varphi\in\mathcal{L}$ and $\tau>0$.\end{prob} The existence of solutions to Problem \ref{problem} indeed provides weak solutions to the following Monge-Amp\`{e}re type equation:
\begin{eqnarray*}
g(\nabla\varphi(y)){\rm det}(\nabla^2\varphi(y))=\tau\cn e^{-\frac{\varphi(y)^2}{2}}e^{-\frac{|y|^2}{2}} 
\end{eqnarray*}
where  $\varphi$ is the unknown function, and   $d\mu=g(y)dy$ with $g$ a smooth function.

\section{A solution to the Euclidean Gaussian Minkowski problem  for convex functions }\label{3.2}

This section aims to solve Problem \ref{problem} when the given measure $\mu$ is   an even measure  and  $\varphi$ is an even function. To this end, let $\mathfrak{M}$ denote the set of all even finite nonzero Borel measures $\mu$ on $\Rn$, such that, $\mu$ is not supported in any lower-dimensional subspaces, and the first moment of $\mu$ is finite, i.e., 
\begin{eqnarray}\label{condition2}
\int_{\Rn}|x|d\mu(x)<\infty.
\end{eqnarray}
Let  ${\rm Supp}(\mu)$ be the support of $\mu$. Denote by ${\rm conv}(E)$ 
the closed convex hull of $E\subset \Rn$. Let  $M_{\mu}$ be the interior of ${\rm conv}({\rm Supp}(\mu))$. Thus, if $\mu\in\mathfrak{M}$, then $ o\in M_{\mu}$.  If $\varphi$ is a $\mu$-integrable convex
function, then $\varphi$ must be finite on $M_{\mu}$.

We consider the following optimization problem: 
\begin{eqnarray}\label{min problem1}
\inf\left\{\int_{\Rn}\varphi(x)d\mu(x):\ \ \varphi \in\mathcal{L}_e^+(\mu) \ \mathrm{and}  \ \G(\varphi^*)=\frac{1}{2} \right\},
\end{eqnarray}
where $\mathcal{L}_e^+(\mu)$ is the class of even, non-negative and $\mu$-integrable functions. Note that, for any $\varphi\in\mathcal{L}_e^+(\mu)$,  $0\leq\varphi^{**}\leq \varphi$ and $\varphi^{***}= \varphi^*$. Then,  
\begin{align*}
\int_{\Rn}\varphi^{**}(x)d\mu(x)\leq \int_{\Rn}\varphi(x)d\mu(x)\quad \mathrm{and}\quad \G(\varphi^*)=\G(\varphi^{***}).
\end{align*}
Consequently,  solving the optimization problem \eqref{min problem1} is equivalent to solving 
\begin{eqnarray}\label{min problem2}
\inf\left\{\int_{\Rn}\varphi(x)d\mu(x):  \ \  \varphi\in \mathcal{L}_e^+(\mu)\cap\mathrm{Conv}(\Rn) \ \ \mathrm{and} \ \ \G(\varphi^*)=\frac{1}{2}\right\}.
\end{eqnarray}

Let $\varphi(x)=a|x|+b$ with $a>0$ and $b>0$. Then, $\varphi^*(x)={\rm I}_{a\ball}^{\infty}(x)-b$. From \eqref{volume1}, one has
\begin{align*}
\G(\varphi^*) =\cn\int_{a\ball}e^{-\frac{|x|^2}{2}}\int_{-b}^{+
\infty}e^{-\frac{s^2}{2}}dsdx  =\cn \bigg(\int_{-b}^{+
\infty}e^{-\frac{s^2}{2}}ds\bigg) \bigg(\int_{a\ball}e^{-\frac{|x|^2}{2}}dx\bigg). 
\end{align*} Clearly, the following identities hold: \begin{align*}
   & \lim _{b\rightarrow +\infty}\int_{-b}^{+
\infty}e^{-\frac{s^2}{2}}ds = \sqrt{2\pi}  \ \ \mathrm{and} \ \ \lim _{b\rightarrow 0 }\int_{-b}^{+
\infty}e^{-\frac{s^2}{2}}ds  =\frac{\sqrt{2\pi}}{2}, \\ 
& \lim_{a\rightarrow +\infty} \int_{a\ball}e^{-\frac{|x|^2}{2}}dx=(2\pi)^{\frac{n}{2}}  \ \ \mathrm{and} \ \  \lim_{a\rightarrow 0} \int_{a\ball}e^{-\frac{|x|^2}{2}}dx=0. 
\end{align*} Consequently, one can find $a_0>0$ and $b_0  >0$ such that $\G(\varphi_0^*)=\frac{1}{2}$ with $\varphi_0(x)=a_0|x|+b_0$.
If $\mu\in\mathfrak{M}$, then the first moment of $\mu$ is finite. Therefore, the optimization problem \eqref{min problem2} is well-defined and 
\begin{align}\label{min}  
\Theta_{\mu} = \inf\left\{\int_{\Rn}\varphi(x)d\mu(x):  \  \varphi\in \mathcal{L}_e^+(\mu)\cap\mathrm{Conv}(\Rn)\  \mathrm{and} \  \G(\varphi^*)=\frac{1}{2}\right\}<\infty.
\end{align}

We shall need the following lemma. 

\begin{lemma}\cite[Lemma 16]{CK15}\label{CK's lemma 16}
Let $\mu$ be a finite Borel measure on $\Rn$. If $x_0\in M_{\mu}$, then there exists $C_{\mu,x_0}>0$ with the following property: for any non-negative, $\mu$-integrable, convex function $\varphi: \Rn\rightarrow [0,\infty]$,
$$
\varphi(x_0)\leq C_{\mu,x_0}\int_{\Rn}\varphi d\mu(x).
$$
\end{lemma} 

The following is another key lemma for our proof. 

\begin{lemma}\cite[Theorem 10.9]{Roc70}\label{Rockafellar's result} Let $C$ be a relatively open convex set, and let $\phi_1, \phi_2,\cdots,$ be a sequence of finite convex functions on $C$. Suppose that the real numbers  $\phi_1(x), \phi_2(x),\cdots,$
	are bounded for each $x\in C$. It is then possible to select a subsequence of  $\phi_1, \phi_2,\cdots,$ which converges to some finite convex function $\phi$ pointwisely on $C$ and uniformly  on  closed  bounded subsets of $C$.
\end{lemma}

The following lemma is similar to  \cite[Lemma 17]{CK15} (see also \cite{FXY,Rot20}), but we make the appropriate modifications according to the need of our main theorem. 

\begin{lemma}\label{existence}
Let $\mu\in \mathfrak{M}$.  If $\varphi_i\in  \mathcal{L}_e^+(\mu)\cap\conv(\Rn)$ and
	\begin{eqnarray}\label{condition1}
	\sup_{i\in \mathbb{N}} \int_{\Rn} \varphi_i(x)d\mu(x)<+\infty.
	\end{eqnarray}
Then, there exists a subsequence $\{\varphi_{i_{j}}\}_{j\in \mathbb{N}}$ of $\{\varphi_i\}_{i\in \mathbb{N}}$ and a function $\varphi\in  \mathcal{L}_e^+(\mu)\cap\conv(\Rn)$ such that
	\begin{align}\label{liminf1}
	\int_{\Rn} \varphi(x) \,d\mu(x) & \leq \liminf_{j\rightarrow\infty} \int_{\Rn}\varphi_{i_j}(x) \,d\mu(x), \\ 
	 \label{liminf2}
	\G (\varphi^*) & \geq \limsup_{j\rightarrow\infty}\G (\varphi_{i_j}^*).
	\end{align}
\end{lemma}

\begin{proof}	  By Lemma \ref{CK's lemma 16} and \eqref{condition1}, Lemma \ref{Rockafellar's result} can be applied to get the existence of a convergence subsequence $\{\varphi_{i_j}\}_{j\in \mathbb{N}}$ of $\{\varphi_i\}_{i\in \mathbb{N}}$, which converges pointwisely to an even non-negative finite convex function $\varphi: M_{\mu} \rightarrow \R$ on $M_{\mu}$ and  converges uniformly on any closed bounded subset of $M_{\mu}$. The function $\varphi: M_{\mu}\rightarrow \R$ can be extended on $\Rn$, still denoted by $\varphi$, by
\begin{eqnarray*}
\varphi(x)=\left\{\begin{array}{lll}
\lim_{\lambda\rightarrow 1^-} \varphi(\lambda x)  & \text{if}&  x\in \partial M_{\mu},  \\  
+\infty  & \text{if}&  x\notin \overline{M_{\mu}}.
\end{array} \right.
\end{eqnarray*} Following the proofs in \cite[Lemma 17]{CK15} (see e.g., \cite[Lemma 5.4]{FYZZ23} and \cite[Lemma 5.8]{HLXZJDG}), one can get inequality \eqref{liminf1} and hence $\varphi\in  \mathcal{L}_e^+(\mu)\cap\conv(\Rn)$.
 
By the continuity of $\varphi$ in $M_{\mu}$,  for any $y \in \Rn$, one gets  $$\varphi^*(y)=\sup_{k\in \mathbb{N}}\big\{\langle x_k,y\rangle-\varphi(x_k)\big\},$$ where $\{x_k\}_{k\in \mathbb{N}}$ is a dense sequence in $M_{\mu}$.  For $j\geq 1$, set $$h_j(y)=\max_{1\leq k \leq j}\big \{\langle x_k,y\rangle-\varphi(x_k)\big\}.$$ Moreover,  $h_j$  is increasing to $\varphi^*$ as $j$ is increasing to $\infty$. It follows from the monotone convergence theorem that
\begin{align*}
	\G(\varphi^*)&=\cn\int_{\Rn}e^{-\frac{|x|^2}{2}} \left(\int_{\varphi^*(x)}^{+
\infty}e^{-\frac{s^2}{2}}ds\right)dx\\
&=\lim_{j\rightarrow\infty}\cn\int_{\Rn}e^{-\frac{|x|^2}{2}} \bigg(\int_{h_j(x)}^{+
\infty}e^{-\frac{s^2}{2}}ds\bigg)dx=\lim_{j\rightarrow\infty}\G(h_j).\end{align*}
 
Let $\varepsilon > 0$. There exists an integer $j_0$ (depending only on $\varepsilon$) satisfying
\begin{align}\label{condition4}
	-\varepsilon\leq \G(h_{j_0})-\G(\varphi^*) \leq\varepsilon.
\end{align} It follows from the fact $\varphi_{i_{j}}\rightarrow \varphi$ pointwisely on $\{x_1,\cdots, x_{j_0}\}$ that  $\varphi^*_{i_j}(x)\geq h_{j_0}(x)-\varepsilon$ holds for all $x\in \Rn$ and for all $j\in \mathbb{N}$ big enough. Together with (\ref{condition4}), one gets that
\begin{align*}
	\G(\varphi^*) &>\G(h_{j_0})-\varepsilon\\& = \cn\int_{\Rn} e^{-\frac{|x|^2}{2}} \int_{h_{j_0}(x)}^{+
\infty}e^{-\frac{s^2}{2}}ds dx-\varepsilon\\& \geq \cn\int_{\Rn} e^{-\frac{|x|^2}{2}} \int_{\varphi^*_{i_j}(x)+\varepsilon}^{+
\infty}e^{-\frac{s^2}{2}}dsdx-\varepsilon
\end{align*}  holds for all $j\in \mathbb{N}$ big enough. Consequently, \begin{align*}
\G(\varphi^*) 	&\geq  \limsup_{j\rightarrow\infty}\cn\int_{\Rn} e^{-\frac{|x|^2}{2}}\bigg(\int_{\varphi^*_{i_j}(x) }^{+
\infty}e^{-\frac{s^2}{2}}ds-\int_{\varphi^*_{i_j}(x)}^{\varphi^*_{i_j}(x)+\varepsilon}e^{-\frac{s^2}{2}}ds\bigg)dx-\varepsilon\\ &\geq    \limsup_{j\rightarrow\infty} \G(\varphi^*_{i_j}) -\limsup_{j\rightarrow\infty}\cn\int_{\Rn}e^{-\frac{|x|^2}{2}}  \int_{\varphi^*_{i_j}(x)}^{\varphi^*_{i_j}(x)+\varepsilon}e^{-\frac{s^2}{2}}ds dx-\varepsilon\\  &\geq    \limsup_{j\rightarrow\infty} \G(\varphi^*_{i_j}) - \big((2\pi)^{-\frac{1}{2}} +1\big) \varepsilon,
\end{align*} where we have used the fact $e^{-\frac{s^2}{2}}\leq 1.$
By letting $\varepsilon\rightarrow 0$, one gets \eqref{liminf2}.
\end{proof}

Now we will deal with the optimization problem \eqref{min problem2}. 

\begin{pro}\label{existence1} For $\mu\in \mathfrak{M}$, there exists a solution $\varphi_0$, which is strictly positive, to the optimization problem \eqref{min problem2}. 
\end{pro}
\begin{proof} Note that the optimization problem \eqref{min problem2} is well-defined and $0\leq \Theta_{\mu}<\infty$ by \eqref{min}.  
We can select a minimizing sequence $\{\varphi_{i}\}_{i\in \mathbb{N}}\in \mathcal{L}_e^+(\mu)\cap\mathrm{Conv}(\Rn)$, such that, for $i\in \mathbb{N}$, 
\begin{align*}
\Theta_{\mu}=\lim_{i\to \infty}\int_{\Rn}\varphi_i(x)d\mu(x) \quad\text{and}\quad \G(\varphi_i^*)=\frac{1}{2}.
\end{align*} In particular,  the condition \eqref{condition1} holds: 
   \begin{eqnarray*}
		\sup_{i\in \mathbb{N}}  \int_{\Rn} \varphi_i(x)\, d\mu(x)<+\infty.
	\end{eqnarray*}  Therefore, Lemma \ref{existence} can be applied to get a subsequence $\{\varphi_{i_{j}}\}_{j\in \mathbb{N}}$ of $\{\varphi_i\}_{i\in \mathbb{N}}$ and  $\varphi_0 \in  \mathcal{L}_e^+(\mu)\cap\conv(\Rn)$ such that
	\begin{align}
		\int_{\Rn} \varphi_0(x) \,d\mu(x) & \leq \liminf_{j\rightarrow\infty} \int_{\Rn} \varphi_{i_j}(x) \,d\mu(x), \nonumber  \\ 
		\G (\varphi_0^*) & \geq \limsup_{j\rightarrow\infty}\G (\varphi_{i_j}^*)=\frac{1}{2}.\label{lim>=1/2}
\end{align}

We now prove that $\varphi_0$ is strictly positive. To this end, assume  that  $\varphi_0(o)=0$. Let $$K_{\varphi_0}=\{x\in \Rn: \varphi_0(x)\leq1\}\ \ \mathrm{and} \ \ r_{\varphi_0}=\min_{v\in \sphere}h_{K_{\varphi_0}}(v).$$ Since $\varphi_0$ is $\mu$-integral, it is finite in a neighborhood of the origin. Together with the convexity of $\varphi_0$, we can obtain $r_{\varphi_0}>0$.
It follows from \eqref{Fenchel conjugate} that \begin{align*}
    \varphi_0^*(y) \geq \sup_{x\in r_{\varphi_0}B_2^n}\left\{\langle x,y\rangle-\varphi_0(x)\right\} \geq \max\big\{-\varphi_0(o),\  r_{\varphi_0}|y|-1\big\}  = \max\big\{0,\  r_{\varphi_0}|y|-1\big\}.  
\end{align*}
More precisely, one has  \begin{align*}\varphi_0^*(y)\geq  \left\{\begin{array}{lll}
0  & \text{if}&  y\in \frac{1}{r_{\varphi_0}}\ball,  \\
r_{\varphi}|y|-1>0  & \text{if}&  y\not\in \frac{1}{r_{\varphi_0}}\ball .
\end{array} \right.
\end{align*}
This further implies that  \begin{align*}
\G(\varphi_0^*)&=\cn \int_{\Rn}e^{-\frac{|x|^2}{2}}\int_{\varphi_0^*(x)}^{+
\infty}e^{-\frac{s^2}{2}}dsdx\\
&<\cn \int_{\Rn}e^{-\frac{|x|^2}{2}}\int_{0}^{+
\infty}e^{-\frac{s^2}{2}}dsdx=\frac{1}{2},
\end{align*} which contradicts to \eqref{lim>=1/2}. Hence, $\varphi_0(o)>0$ and then $\varphi_0$ is strictly
positive as $\varphi_0 \in  \mathcal{L}_e^+(\mu)\cap\conv(\Rn)$. 

\vskip 2mm Next, we prove that  $\G (\varphi_0^*) = \frac{1}{2}$, again by the argument of contradiction. That is, we assume  $\G (\varphi_0^*) > \frac{1}{2}$. For any   $\tau\geq 0$,  let $\varphi_{\tau}=\max\{0, \varphi_0-\tau\}$ be a nonnegative, even convex function.  Clearly, $\varphi_0-\tau \leq \varphi_{\tau}\leq \varphi_{0}$ as $\varphi$ is strictly positive, and hence  by \eqref{moto} and \eqref{star}, one has $\varphi_0^* \leq \varphi_{\tau}^* \leq \varphi_0^*+\tau$ for any $\tau\geq 0$.  Also note that $\varphi_{\tau}$ is decreasing  and hence $\varphi_{\tau}^*$ is increasing  on $\tau \geq 0$. Thus, $\gamma_{n+1}(\varphi_{\tau}^*)$ is decreasing on $\tau>0.$ Another useful fact is that $\dom \varphi_{\tau}=\dom  \varphi_0$ for any $\tau>0.$

On the one hand, as $\varphi_0$ is strictly positive, for any $0<\tau<\varphi_0(o)$, $\varphi_{\tau}=\varphi_0-\tau$ and then by \eqref{star}, $\varphi_{\tau}^*=\varphi_0^*+\tau.$ This further gives, for any $0<\tau<\varphi_0(o)$,  \begin{align*} \G(\varphi_{0}^*)&\geq 
\G(\varphi_{\tau}^*) \nonumber\\ &  =\cn \int_{\Rn}e^{-\frac{|x|^2}{2}}\int_{\varphi_0^*(x)+ \tau }^{+
\infty}e^{-\frac{s^2}{2}}dsdx \nonumber \\
&=\cn \int_{\Rn}e^{-\frac{|x|^2}{2}}\int_{\varphi_0^*(x)}^{+
\infty}e^{-\frac{s^2}{2}}dsdx -\cn \int_{\Rn}e^{-\frac{|x|^2}{2}}\int^{\varphi_0^*(x)+\tau}_{\varphi_0^*(x)}e^{-\frac{s^2}{2}}dsdx \nonumber \\ 
&\geq \G(\varphi_0^*) -\frac{\tau}{\sqrt{2\pi}}, 
\end{align*}  where again we have used $0< e^{-\frac{s^2}{2}}\leq 1$ for all $s\in \R.$  Clearly, 
\begin{align}
    \label{lim>1/2} \lim_{ \tau\rightarrow  0^+ } \G(\varphi_{\tau}^*)= \G(\varphi_0^*)>\frac{1}{2}. 
\end{align}

On the other hand, as $\varphi_0 \in  \mathcal{L}_e^+(\mu)\cap\conv(\Rn)$, then $\varphi_0$ must be finite on $M_{\mu}$. Note that  $o\in M_{\mu}$ for $\mu\in \mathfrak{M}$. Let $r_0>0$ be such that $r_0B^n_2\subset M_{\mu}$, and $$\tau_0=\max\big\{\varphi_0(x): x\in r_0B^n_2\big\}>0.$$  It can be checked that $\varphi_{\tau_0}\leq {\rm I}^{\infty}_{r_0B^n_2}$ and hence, $$\varphi_{\tau_0}^*\geq ({\rm I}^{\infty}_{r_0B^n_2})^*=h_{r_0B^n_2}=r_0|\cdot|.$$  This further implies that \begin{align}
    \G\big(\varphi_{\tau_0}^*\big) &=\cn \int_{\Rn}e^{-\frac{|x|^2}{2}}\int_{\varphi_{\tau_0}^*(x)}^{+
\infty}e^{-\frac{s^2}{2}}dsdx  \nonumber \\&\leq  \cn \int_{\Rn}e^{-\frac{|x|^2}{2}}\int_{r_0|x|}^{+
\infty}e^{-\frac{s^2}{2}}dsdx <\frac{1}{2}.  \label{H-limit=0}
\end{align} 
 
We now prove that there exists a $\tau_1\in (0, \tau_0)$ such that \begin{align*} \G\big(\varphi_{\tau_1}^*\big) =\frac{1}{2}.
\end{align*} To this end, we need to show the continuity of $ \G\big(\varphi_{\tau}^*\big)$ on $\tau>0.$ Let $\tau>0$ be any given number and $0<\delta_0<\frac{\tau}{2}$. For any $t$ such that $|t-\tau|<\delta_0$ (i.e., $0<\tau-\delta_0<t<\tau+\delta_0)$), one has, for all $x\in \dom \varphi_0$,    
$$|\varphi_t-\varphi_{\tau}| =\big|\max\{0, \varphi_0 -t\}-\max\{0, \varphi_0 -\tau\}\big|\leq |t-\tau|. $$ This further yields that, for any $x\in \Rn,$ $$\varphi_{\tau}(x)-|t-\tau|\leq \varphi_{t} (x)\leq \varphi_{\tau}(x)+|t-\tau|. $$ It follows from \eqref{moto} and  \eqref{star} that 
\begin{align*}
    \varphi_{\tau}^*-|t-\tau|\leq \varphi_{t} ^* \leq \varphi_{\tau}^*+|t-\tau|.
\end{align*} Together with formula \eqref{volume1}, one has  \begin{align*}
    \big| \G(\varphi_{\tau}^*)-\G(\varphi_{t}^*) \big| &= \cn \bigg|   \int_{\Rn }e^{-\frac{|x|^2}{2}}\int_{\varphi_{\tau}^*(x)}^ {\varphi_{t}^*(x)}e^{-\frac{s^2}{2}}dsdx\bigg| \nonumber \\  &\leq \cn   \int_{\Rn }e^{-\frac{|x|^2}{2}} \bigg|  \int_{\varphi_{\tau}^*(x)}^ {\varphi_{t}^*(x)}e^{-\frac{s^2}{2}}ds \bigg|  dx   \nonumber  \\  &\leq \cn  |t-\tau|  \int_{\Rn }e^{-\frac{|x|^2}{2}}   dx\nonumber  \\  & =   \frac{|t-\tau|}{\sqrt{2\pi}},
\end{align*} where we have used $e^{-\frac{s^2}{2}}\leq 1$ in the last inequality. This immediately yields the continuity  of $ \G\big(\varphi_{\tau}^*\big)$ on $\tau>0.$  Together with \eqref{lim>1/2} and \eqref{H-limit=0},  one can find  $\tau_1>0$, such that  $ \G(\varphi_{\tau_1}^*)= \frac{1}{2}$ and $ \varphi_0 - \varphi_{\tau_1}>0$.   The latter one yields that $$\int_{\Rn} \varphi_0(x) \,d\mu(x)>\int_{\Rn} \varphi_{\tau_1}(x) \,d\mu(x),$$ which contradicts to the minimality of $\int_{\Rn}\varphi_0\,d\mu$ (in view of $\G(\varphi_{\tau_1}^*)= \frac{1}{2}$). Therefore, $\G (\varphi_0^*) = \frac{1}{2},$ and then $\varphi_0$ solves the optimization problem \eqref{min problem2} (and hence, \eqref{min problem1}). 
\end{proof}

In the last part of this section, we will prove that, if  the convex function $\varphi_0$ solves  the optimization problem \eqref{min problem2}, then $\varphi_0$ is a solution to the Euclidean Gaussian Minkowski problem of convex functions (i.e., Problem \ref{problem}). 
The following result is needed.
\begin{lemma}\label{dominate}
Let $\varphi:\Rn\to\R\cup\{+\infty\}$ be a lower semi-continuous function  with $\varphi(o)<+\infty$. Assume that $g:\Rn\to\R$ is bounded and continuous. Then
\begin{eqnarray*}
\frac{d}{dt}\Big|_{t=0} \G((\varphi+tg)^*)
=\cn \int_{\Rn}g(\nabla\varphi^*(x))e^{-\frac{\varphi^*(x)^{2}}{2}}e^{-\frac{|x|^2}{2}}dx.
\end{eqnarray*}
\end{lemma}
\begin{proof}
Applying Lemma \ref{Berman's formula} to $g$ and $-g$ at any point $x\in \Rn$ in which $\varphi^*$ is differentiable, one gets 
 
\begin{align}\label{difflem}
\frac{d}{dt}\Big |_{t=0}(\varphi+tg)^*(x)=-g(\nabla\varphi^*(x)).
\end{align}

Assume that $|g|\leq M$ for some $M>0$.
Then, \begin{align*}
\varphi-|t|M \leq\varphi+tg\leq\varphi+|t|M.
\end{align*}
From \eqref{moto} and  \eqref{star}, one has
\begin{align}\label{varphicontinuity}
\varphi^*-|t|M \leq(\varphi+tg)^*\leq\varphi^*+|t|M.
\end{align}
As $e^{-\frac{s^2}{2}}\leq 1$,   for any $x\in\Rn$, one has 
\begin{align*}
\frac{1}{t}\left| \int_{(\varphi+tg)^*(x)}^{+
\infty}e^{-\frac{s^2}{2}}ds- \int_{\varphi^*(x)}^{+
\infty}e^{-\frac{s^2}{2}}ds\right| = \frac{1}{t}\left| \int_{(\varphi+tg)^*(x)}^ {\varphi^*(x)} e^{-\frac{s^2}{2}}ds\right| \leq M.
\end{align*}
Together with \eqref{difflem}, the dominated convergence theorem deduces that
\begin{align*}
\frac{d}{dt}\Big|_{t=0} \G((\varphi+tg)^*)&= \lim_{t\to0}\frac{\G( (\varphi+tg)^*)-\G(  \varphi^*)}{t}\nonumber\\
&=\cn \lim_{t\to0}\int_{\Rn}
\frac{1}{t}\bigg(  \int_{(\varphi+tg)^*(x)}^{+
\infty}e^{-\frac{s^2}{2}}ds- \int_{\varphi^*(x)}^{+
\infty}e^{-\frac{s^2}{2}}ds \bigg)  e^{-\frac{|x|^2}{2}}dx\nonumber\\&=\cn \int_{\Rn}
\lim_{t\to0} \frac{1}{t}\bigg(  \int_{(\varphi+tg)^*(x)}^ {\varphi^*(x)} e^{-\frac{s^2}{2}}ds \bigg)  e^{-\frac{|x|^2}{2}}dx\nonumber\\
&=\cn \int_{\Rn}g(\nabla\varphi^*(x))e^{-\frac{\varphi^*(x)^{2}}{2}}e^{-\frac{|x|^2}{2}}dx.
\end{align*} This completes the proof. \end{proof}

We now prove our main result, the existence of solution to Problem \ref{problem}. 

\begin{thm}\label{main-existence} Let $\mu\in \mathfrak{M}$. Then  there exists  $\varphi\in\mathcal{L}$ such that
\begin{align}\label{final-result}
 d\mu  =\frac{|\mu|}{\mu_{\gamma_n}(\varphi, \Rn)}d \mu_{\gamma_n}(\varphi, \cdot),
 \end{align}  where $|\mu|$ and $\mu_{\gamma_n}(\varphi, \Rn)$ are real numbers given by 
 $$\mu_{\gamma_n}(\varphi, \Rn)=\int_{\Rn} d \mu_{\gamma_n}(\varphi, x) \ \ \mathrm{and} \ \ |\mu|=\int_{\Rn}\,d\mu. $$
\end{thm}
\begin{proof}
According to Proposition \ref{existence1}, there exists $\varphi_0\in\mathcal{L}_e^+(\mu)\cap\conv(\Rn)$ solving the optimization problem \eqref{min problem2}. Moreover $\varphi_0>0$.

Let $g:\Rn\rightarrow \R$ be an even compactly supported continuous function. Then, $g$ is bounded on $\Rn$, i.e., $|g|<M$ for some $M$. For $t_1, t_2\in \R$, let \begin{align}
  \varphi_{t_1,t_2}(x)=\varphi_0(x)+t_1g(x)+t_2. \label{two-parameter-varphi}  
\end{align} As $\varphi_0>0$, for sufficiently small $t_0,t_0'>0$, $\varphi_{t_1,t_2}(x)\in\mathcal{L}_e^+(\mu)$ for $t_1\in[-t_0,t_0]$ and $t_2\in[-t_0',t_0']$. Consequently, for sufficiently small $t$,  \begin{align*}
    \varphi_{t_1+t,t_2}(x) =\varphi_{t_1,t_2}(x)+tg(x) \ \ \mathrm{and}\ \  \varphi_{t_1,t_2+t}(x)=\varphi_{t_1,t_2}(x)+t,
\end{align*} which are both in $\mathcal{L}_e^+(\mu)$.  
Applying Lemma \ref{dominate} (to $\varphi=\varphi_{t_1,t_2}$), one gets 
\begin{align}\label{partial1}
\frac{\partial}{\partial t_1}\G(\varphi_{t_1,t_2}^*)&= \lim_{t\to0}\frac{\G( \varphi_{t_1+t,t_2}^*)-\G(  \varphi_{t_1,t_2}^*)}{t}\nonumber\\
&=\lim_{t\to0}\frac{\G( (\varphi_{t_1,t_2}+tg)^*)-\G(  \varphi_{t_1,t_2}^*)}{t}\nonumber\\
&=\cn \int_{\Rn}g(\nabla\varphi_{t_1,t_2}^*(x))e^{-\frac{\left(\varphi_{t_1,t_2}^*(x)\right)^{2}}{2}}e^{-\frac{|x|^2}{2}}dx.
\end{align} 
Similarly,   Lemma  \ref{dominate} implies 
\begin{align}\label{partial2}
\frac{\partial}{\partial t_2}\G(\varphi_{t_1,t_2}^*)&= \lim_{t\to0}\frac{\G( \varphi_{t_1,t_2+t}^*)-\G(  \varphi_{t_1,t_2}^*)}{t}\nonumber\\
&=\lim_{t\to0}\frac{\G( (\varphi_{t_1,t_2}+t)^*)-\G(  \varphi_{t_1,t_2}^*)}{t}\nonumber\\
&=\cn \int_{\Rn}e^{-\frac{\left(\varphi_{t_1,t_2}^*(x)\right)^{2}}{2}}e^{-\frac{|x|^2}{2}}dx\nonumber\\
&=\mu_{\gamma_n}(\varphi_{t_1,t_2}^*, \Rn).
\end{align} 

Now we claim that both $\frac{\partial}{\partial t_1}\G(\varphi_{t_1,t_2}^*)$ and $\frac{\partial}{\partial t_2}\G(\varphi_{t_1,t_2}^*)$ are continuous on $(t_1,t_2)\in S_0$ with $S_0=[-t_0,t_0]\times [-t_0',t_0'].$ Let $(t_1, t_2)\in S_0$, and let $\{r_i\}_{i\in \N}$ and $\{s_i\}_{i\in \N}$ be sequences convergent to $0$  such that $(t_1+r_i, t_2+s_i)\in S_0$ for all $i\in \N$.  Following the proof for  \eqref{varphicontinuity}, and by $|g|\leq M$ on $\Rn$, one has, for any $i\in \N,$ 
\begin{align*}
\varphi_{t_1,t_2}^*-|r_i|M-|s_i|\leq\varphi^*_{t_1+r_i,t_2+s_i}=(\varphi_{t_1,t_2}+r_ig+s_i)^*\leq\varphi_{t_1,t_2}^*+|r_i|M+|s_i|.
\end{align*}
This further implies that \begin{align}\label{t1t2continuity}\lim_{i  \rightarrow \infty}\varphi^*_{t_1+r_i,t_2+s_i}=\varphi^*_{t_1, t_2}.\end{align}   Moreover, for $i\in \N$, one has $$D_{\varphi_{t_1+r_i, t_2+s_i}^*}=D_{\varphi_{t_1, t_2}^*}=D_{\varphi_0^*}.$$ 
It follows from \cite[Theorems 24.5]{Roc70} that $\nabla\varphi_{t_1+r_i,t_2+s_i}^*(x)$ converges pointwisely to $\nabla\varphi_{t_1,t_2}^*(x)$ at those $x$ where $\nabla\varphi_{t_1+r_i,t_2+s_i}^*(x)$ for $i\in \N$ are all differentialble. Note that, for each $i\in \N$, $ \varphi_{t_1+r_i,t_2+s_i}^*$ is differentiable almost everywhere in  $\rm {int}(D_{\varphi_0^*})$, and hence $\nabla\varphi_{t_1+r_i,t_2+s_i}^*(x)$ converges pointwisely to $\nabla\varphi_{t_1,t_2}^*(x)$ almost everywhere in $\rm {int}(D_{\varphi_0^*})$.
 As $g$ is an even compactly supported continuous function,   by \eqref{t1t2continuity}, one has, for almost any $x\in\Rn$,
\begin{align}\label{limit-pointwise-1}
    \lim_{i\to \infty}g\big(\nabla\varphi_{t_1+r_i,t_2+s_i}^*(x)\big)e^{-\frac{\left(\varphi_{t_1+r_i, t_2+s_i}^*(x)\right)^{2}}{2}}=g\big(\nabla\varphi_{t_1,t_2}^*(x)\big)e^{-\frac{\left(\varphi_{t_1,t_2}^*(x)\right)^{2}}{2}}. 
\end{align} By $|g|\leq M$ on $\Rn$, one has, for almost any $x\in\Rn$ and for all $i\in \N$, $$ \Big| g(\nabla\varphi_{t_1+r_i,t_2+s_i}^*(x))\Big| e^{-\frac{\left(\varphi_{t_1+r_i,t_2+s_i}^*(x)\right)^{2}}{2}}e^{-\frac{|x|^2}{2}}\leq Me^{-\frac{|x|^2}{2}}.$$ It follows from \eqref{partial1}, \eqref{limit-pointwise-1} and the dominated convergence theorem that 
\begin{align*} \frac{\partial}{\partial t_1}\G(\varphi_{t_1,t_2}^*)&=\cn    \int_{\Rn}\! g(\nabla\varphi_{t_1,t_2}^*(x))e^{-\frac{\left(\varphi_{t_1,t_2}^*(x)\right)^{2}}{2}}e^{-\frac{|x|^2}{2}}dx\\ &=\cn \lim_{i\to \infty}\int_{\Rn} g(\nabla\varphi_{t_1+r_i,t_2+s_i}^*(x))e^{-\frac{\left(\varphi_{t_1+r_i,t_2+s_i}^*(x)\right)^{2}}{2}}e^{-\frac{|x|^2}{2}}dx \\ &=  \lim_{i\to \infty}  \frac{\partial}{\partial t_1}\G(\varphi_{t_1+r_i,t_2+s_i}^*).\end{align*}  As the sequences $\{r_i\}_{i\in \N}$ and $\{s_i\}_{i\in \N}$ are arbitrary, one gets that $\frac{\partial}{\partial t_1}\G(\varphi_{t_1,t_2}^*)$ is continuous on $(t_1,t_2)\in S_0$. Similarly, for each $i\in \N$ and $x\in \Rn$, it holds that
$$e^{-\frac{\left(\varphi_{t_1+r_i,t_2+s_i}^*(x)\right)^{2}}{2}}e^{-\frac{|x|^2}{2}}\leq e^{-\frac{|x|^2}{2}}.$$ Again, due to \eqref{partial2}, \eqref{t1t2continuity} and the dominated convergence theorem, one gets 
\begin{align*} 
\frac{\partial}{\partial t_2}\G(\varphi_{t_1,t_2}^*) 
&=\cn \int_{\Rn}e^{-\frac{\left(\varphi_{t_1,t_2}^*(x)\right)^{2}}{2}}e^{-\frac{|x|^2}{2}}dx\nonumber\\ &=\cn \lim_{i\rightarrow \infty} \int_{\Rn}e^{-\frac{\left(\varphi_{t_1+r_i,t_2+s_i}^*(x)\right)^{2}}{2}}e^{-\frac{|x|^2}{2}}dx\nonumber\\ &=\lim_{i\rightarrow \infty} \frac{\partial}{\partial t_2}\G(\varphi_{t_1+r_i,t_2+s_i}^*). 
\end{align*} As the sequences $\{r_i\}_{i\in \N}$ and $\{s_i\}_{i\in \N}$ are arbitrary, one gets that $\frac{\partial}{\partial t_2}\G(\varphi_{t_1,t_2}^*)$ is continuous on $(t_1,t_2)\in S_0$. 

On the other hand, one notices that, for any $(t_1, t_2)\in S_0,$
$$\frac{\partial}{\partial t_2} \G(\varphi_{t_1,t_2}^*)=
\cn \int_{\Rn}e^{-\frac{\left(\varphi_{t_1, t_2}^*(x)\right)^{2}}{2}}e^{-\frac{|x|^2}{2}}dx > 0.$$ This is an easy consequence from $\G(\varphi^*_0)=\frac{1}{2}$, yielding that $D_{\varphi_0^*}$ (and hence $D_{\varphi_{t_1, t_2}^*}$) has positive Lebesgure measure.  
These allow us to use the Lagrange multiplier method to the optimization problem \eqref{min problem2}. To this end, for $t_1, t_2, \lambda\in \R$, let $$\Psi(t_1,t_2,\lambda)=\int_{\Rn} \varphi_{t_1,t_2}(x) \,d\mu(x)+\lambda\bigg(\frac{1}{2}-\G(\varphi_{t_1,t_2}^*)\bigg).$$ 
As $\varphi_0$ solves the optimization problem \eqref{min problem2},   the Lagrange multiplier method implies  that
\begin{align*}
 \frac{\partial}{\partial t_1}\Big |_{t_1=t_2=0} \Psi(t_1,t_2,\lambda)=0 \ \ \mathrm{and}\ \  \frac{\partial}{\partial t_2}\Big |_{t_1=t_2=0} \Psi(t_1,t_2,\lambda)=0.    
\end{align*} Consequently, the following equations hold:  
\begin{align}\label{domian1}
\frac{\partial}{\partial t_1}\Big |_{t_1=t_2=0} \left( \int_{\Rn}\varphi_{t_1,t_2}(x)d\mu(x)\right)
=\lambda\frac{\partial}{\partial t_1}\Big |_{t_1=t_2=0}\G(\varphi_{t_1,t_2}^*),\\ 
 \label{domian2}
\frac{\partial}{\partial t_2}\Big |_{t_1=t_2=0} \left( \int_{\Rn}\varphi_{t_1,t_2}(x)d\mu(x)\right) = \lambda\frac{\partial}{\partial t_2}\Big |_{t_1=t_2=0}  \G(\varphi_{t_1,t_2}^*).
\end{align} Due to \eqref{two-parameter-varphi}, it is easily checked that  $$\int_{\Rn}\varphi_{t_1,t_2}(x)d\mu(x)=\int_{\Rn}\varphi_0(x)d\mu(x)+t_1\int_{\Rn}g(x)d\mu(x)+t_2\int_{\Rn}d\mu(x).$$ Thus, the following identities can be obtained: 
\begin{align} 
\frac{\partial}{\partial t_1}\Big |_{t_1=t_2=0}\left( \int_{\Rn}\varphi_{t_1,t_2}(x)d\mu(x)\right)
&=\int_{\Rn}g(x)d\mu(x), \nonumber \\ 
\frac{\partial}{\partial t_2}\Big |_{t_1=t_2=0}\left( \int_{\Rn}\varphi_{t_1,t_2}(x)d\mu(x)\right)&={|\mu|}. \label{derivative=constant}
\end{align} Together with  \eqref{partial1} and \eqref{domian1}, one can conclude that, for any even compactly supported continuous function $g$,  \begin{align}  
\int_{\Rn}g(x)d\mu(x)   =\lambda\cn \int_{\Rn}g(\nabla\varphi_{0}^*(x))e^{-\frac{\varphi_{0}^*(x)^{2}}{2}}e^{-\frac{|x|^2}{2}}dx =\lambda\int_{\Rn}g(x)d\mu_{\gamma_n}(\varphi_0^*, x). \label{equation-g}  
\end{align} Similarly, by \eqref{partial2},  \eqref{domian2} and \eqref{derivative=constant}, one gets  $|\mu|  =\lambda\mu_{\gamma_n}(\varphi_0^*, \Rn).$ Thus, $\lambda\in \R$ is a fixed constant independent of $g$, namely,  
 $$\lambda=\frac{|\mu|}{\mu_{\gamma_n}(\varphi_0^*, \Rn)}.$$ This, together with \eqref{equation-g}, yields that \begin{align} 
 d\mu  =\frac{|\mu|}{\mu_{\gamma_n}(\varphi_0^*, \Rn)}d \mu_{\gamma_n}(\varphi_0^*, \cdot). \label{solution-varphi-0}
 \end{align}

Note that $\varphi_0$ is finite in a neighborhood of the origin, and thus, $o\in \rm int (dom~\varphi_0)$. It follows from  \cite[Theorem 11.8 (c)]{Roc98} that $\varphi_0^*\in\mathcal{L}$. If we let $\varphi=\varphi_0^*$, then $\varphi\in\mathcal{L}$ is a proper, even and lower semi-continuous convex function. In particular,  \eqref{solution-varphi-0} can be written by 
\begin{align*} 
 d\mu  =\frac{|\mu|}{\mu_{\gamma_n}(\varphi, \Rn)}d \mu_{\gamma_n}(\varphi, \cdot),
 \end{align*} which is the desired formula \eqref{final-result}. This completes the proof.  
\end{proof}

{\bf Acknowledgment.} 
The research of XL has been supported by the Science and Technology Research Program of Chongqing Municipal Education Commission (No. KJQN202300557) and the  Research Foundation of Chongqing Normal University (No. 20XLB012). The research of DY has been supported by a NSERC grant.

\end{document}